\documentclass{article}

\usepackage{graphicx}
\usepackage{color}
\usepackage{amssymb,amsmath}
\usepackage{mathrsfs}
\usepackage{bbm}
\usepackage[text={144mm,220mm},centering]{geometry}

\usepackage[numbers]{natbib}
\usepackage{natbibspacing}
\renewcommand{\bibfont}{\small}

\definecolor{myurlcolor}{rgb}{0,0,0.5}
\usepackage{hyperref}
\hypersetup{colorlinks,
linkcolor=black,
citecolor=black,
urlcolor=myurlcolor}

\newcommand{\hyph}{\mbox{-}}
\newcommand{\ucontents}[2]{\addcontentsline{toc}{#1}{\numberline{}{#2}}}
\newcommand{\cat}[1]{\mathscr{#1}}
\newcommand{\fcat}[1]{\mathbf{#1}}
\newcommand{\slsh}{/\linebreak[0]}
\newcommand{\dblslsh}{//\linebreak[0]}
\newcommand{\dt}{.\linebreak[0]}
\newcommand{\such}{:}
\newcommand{\without}{\setminus}
\newcommand{\epsln}{\varepsilon}
\newcommand{\Hom}{\mathrm{Hom}}
\newcommand{\Cat}{\fcat{Cat}}
\newcommand{\Set}{\fcat{Set}}
\newcommand{\integers}{\mathbb{Z}}
\newcommand{\goesto}{\mapsto}

\newcommand{\oppairu}{\rightleftarrows}
\newcommand{\reals}{\mathbb{R}}
\newcommand{\rationals}{\mathbb{Q}}
\newcommand{\complexes}{\mathbb{C}}
\newcommand{\latin}[1]{\textit{#1}}	
\newcommand{\demph}[1]{\textbf{\textup{#1}}}
\newcommand{\done}{\hfill\ensuremath{\Box}}
\newenvironment{prooflike}[1]{\begin{trivlist}\item\textbf{#1}\ }
{\end{trivlist}}
\newenvironment{proof}{\begin{prooflike}{Proof}}{\end{prooflike}}
\newcommand{\scat}[1]{\mathbf{#1}}

\newcommand{\iso}{\cong}
\newcommand{\nat}{\mathbb{N}}	
\newcommand{\eqv}{\simeq}
\DeclareMathOperator{\ob}{ob}		
\newcommand{\pr}{\mathrm{pr}}
\newcommand{\of}{\circ}
\newcommand{\sub}{\subseteq}
\newcommand{\cell}[4]{\put(#1,#2){\makebox(0,0)[#3]{\ensuremath{#4}}}}
\newcommand{\zmark}{\scriptstyle{\bullet}}
\newcommand{\tr}{\texttt{t}}
\newcommand{\fa}{\texttt{f}}
\newcommand{\Two}{\mathbf{2}}
\newcommand{\slogan}[1]%
{\begin{center}\textit{#1}\end{center}}
\newcommand{\restr}[1]{\vert_{#1}}
\newcommand{\Setf}{\fcat{FinSet}}
\newcommand{\Vectf}{\fcat{FDVect}}
\newcommand{\VCat}{\cat{V}\!\hyph\Cat}
\newcommand{\VCatf}{\cat{V}\!\hyph\fcat{cat}}
\newcommand{\MS}{\fcat{MS}}
\newcommand{\Vect}{\fcat{Vect}}
\newcommand{\mg}[1]{|#1|}
\newcommand{\mgfml}[1]{[ #1 ]}
\newcommand{\nm}[1]{\Vert #1 \Vert}
\newcommand{\nrm}[2]{\Vert #1 \Vert_{#2}}
\newcommand{\ft}[1]{\hat{#1}}
\newcommand{\cvx}[1]{\mathscr{K}_{#1}}
\newcommand{\cvxa}[1]{\mathscr{K}(#1)}
\newcommand{\Val}[1]{\mathrm{Val}_{#1}}
\newcommand{\Vala}[1]{\mathrm{Val}(#1)}
\DeclareMathOperator{\Vol}{Vol}
\newcommand{\Gr}[2]{\mathrm{Gr}_{#1,#2}}
\newcommand{\onedec}[1]{#1'}
\newcommand{\Grone}[2]{\onedec{\Gr{#1}{#2}}}
\newcommand{\Vone}{\onedec{V}}
\newcommand{\Valone}[1]{\onedec{\Val{#1}}}
\newcommand{\cvxone}[1]{\onedec{\mathscr{K}_{#1}}}
\newcommand{\Gram}[2]{\nu_{#1,#2}}
\newcommand{\mdim}{\dim}
\DeclareMathOperator{\sech}{sech}
\newcommand{\tra}{*}
\newcommand{\ip}[2]{\langle #1, #2 \rangle}
\newcommand{\dee}{\,d}
\DeclareMathOperator{\adj}{adj}
\DeclareMathOperator{\Nil}{Nil}
\newcommand{\Unit}{\mathbbm{1}}

\newcommand{\incl}{\hookrightarrow}
\newcommand{\partialmap}{\dashrightarrow}
\newcommand{\GL}[2]{\mathrm{GL}(#1)}
\newcommand{\Orthog}[1]{\mathrm{O}(#1)}
\newcommand{\Field}[1]{\mathbb{F}_{#1}}
\newcommand{\cln}{\colon}
\newcommand{\card}[1]{\# #1}
\newcommand{\minev}{\lambda_{\text{min}}}
\newcommand{\UM}[2]{U_{#1}^{#2}}
\newcommand{\Sch}{\mathscr{S}}
\newcommand{\Poset}{\fcat{Poset}}
\newcommand{\cuboidclass}[1]{\mathscr{C}_{#1}}
\newcommand{\fnbound}[1]{c(#1)}
\newcommand{\spacebound}[1]{\langle #1\rangle}
\newcommand{\WUM}{\Omega}
\DeclareMathOperator{\diam}{diam}
\DeclareMathOperator{\Sing}{Sing}
\newcommand{\pv}[1]{\mathbf{#1}}
\newcommand{\Dmax}[1]{D_\text{max}(#1)}
\renewcommand{\implies}{\Rightarrow}
\newcommand{\mrk}[1]{#1}
\newcommand{\codi}[3]{#1 +_{#2} #3}

\newtheorem{thm}{Theorem}[subsection]
\newtheorem{propn}[thm]{Proposition}
\newtheorem{lemma}[thm]{Lemma}
\newtheorem{cor}[thm]{Corollary}
\newtheorem{conj}[thm]{Conjecture}
\newtheorem{predefn}[thm]{Definition}
\newenvironment{defn}{\begin{predefn}\upshape}{\end{predefn}}
\newtheorem{preexample}[thm]{Example}
\newenvironment{example}{\begin{preexample}\upshape}{\end{preexample}}
\newtheorem{preexamples}[thm]{Examples}
\newenvironment{examples}{\begin{preexamples}\upshape}{\end{preexamples}}

\author{Tom Leinster%
\thanks{School of Mathematics and Statistics, University of Glasgow, Glasgow
G12 8QW, UK; Tom.Leinster@glasgow.ac.uk.  
Supported by an EPSRC Advanced Research Fellowship.}}
\title{\vspace*{-2cm}The magnitude of metric spaces}
\date{}

\begin{document}

\sloppy
\maketitle

\begin{abstract}
Magnitude is a real-valued invariant of metric spaces, analogous to the Euler
characteristic of topological spaces and the cardinality of sets.  The
definition of magnitude is a special case of a general categorical definition
that clarifies the analogies between cardinality-like invariants in
mathematics.  Although this motivation is a world away from geometric measure,
magnitude, when applied to subsets of $\reals^n$, turns out to be intimately
related to invariants such as volume, surface area, perimeter and dimension.
We describe several aspects of this relationship, providing evidence for a
conjecture (first stated in~\cite{AMSES}) that magnitude subsumes all the most
important invariants of classical integral geometry.
\end{abstract}

\tableofcontents\ 

\section*{Introduction}
\ucontents{section}{Introduction}

Many mathematical objects carry a canonical notion of size.  Sets have
cardinality, vector spaces have dimension, topological spaces have Euler
characteristic, and probability spaces have entropy.  This work adds a new
item to the list: metric spaces have magnitude.

Already, several cardinality-like invariants are tied together by the notion
of the Euler characteristic of a category \cite{ECC,ECCSDS}.  This is a
rational-valued invariant of finite categories.  A network of theorems
describes the close relationships between this invariant and established
cardinality-like invariants, including the cardinality of sets and of
groupoids~\cite{BaDoFSF}, the Euler characteristic of topological spaces and
of posets, and even the Euler characteristic of orbifolds.  (That Euler
characteristic deserves to be considered an analogue of cardinality was first
made clear by Schanuel~\cite{SchaWLP,SchaNSE}.)  These results attest that for
categories, Euler characteristic is the fundamental notion of size.

Here we go further.  Categories are a special case of the more general concept
of enriched category.  Much of ordinary category theory generalizes to the
enriched setting, and this is true, in particular, of the Euler characteristic
of categories.  Rebaptizing Euler characteristic as `magnitude' to avoid
a potential ambiguity, this gives a canonical definition of the magnitude of
an enriched category.

Metric spaces, as well as categories, are examples of enriched
categories:
\begin{center}
(categories) $\subset$ (enriched categories) $\supset$ (metric spaces)
\end{center}
\cite{LawMSG,LawTCS}.
The analogy between categories and metric spaces can be understood
immediately.  A category has objects; a metric space has points.  For any two
objects there is a set (the maps between them); for any two points there is a
real number (the distance between them).  For any three objects there is an
operation of composition; for any three points there is a triangle inequality.

Having generalized the definition of magnitude (or Euler characteristic) from
ordinary to enriched categories, we specialize it to metric spaces.  This
gives our invariant.  The fundamental role of the Euler characteristic of
categories strongly suggests that the magnitude of metric spaces should play a
fundamental role too.  Our faith is rewarded by a series of theorems showing
that magnitude is intimately related to the classical invariants of integral
geometry: dimension, perimeter, surface area, volume, \ldots.  This is despite
the fact that no concept of measure or integration goes into the definition of
magnitude; they arise spontaneously from the general categorical definition. 

While the author's motivation was category-theoretic, magnitude had already
arisen in work on the quantification of biodiversity.  In 1994, Solow and
Polasky~\cite{SoPo} carried out a probabilistic analysis of the benefits of
high diversity, and isolated a particular quantity that they called the
`effective number of species'.  It is the same as our magnitude.  As it
transpires, this is no coincidence: under suitable
circumstances~\cite{METAMB}, magnitude can be interpreted as maximum
diversity, a cousin to maximum entropy.

We start by defining the magnitude of an enriched category
(Section~\ref{sec:enr}).  This puts the notion of the magnitude of a metric
space into a wide mathematical context, showing how analogous theories can be
built in parts of mathematics far away from metric geometry.  The reader
interested only in geometry can, however, avoid these general considerations
without logical harm.  Such a reader can begin at Section~\ref{sec:fin}.

A topological space is not guaranteed to have a well-defined Euler
characteristic unless it satisfies some finiteness condition.  Similarly, the
magnitude of an enriched category is defined under an assumption of
finiteness; specializing to metric spaces, the definition of magnitude
is just for \emph{finite} spaces (Section~\ref{sec:fin}).  The magnitude of a
finite metric space can be thought of as the `effective number of points'.  It
deserves study partly because of its intrinsic interest, partly because of its
applications to the measurement of diversity, and partly because it is used in
the theory of magnitude of \emph{infinite} metric spaces.

While categorical arguments do not (yet) furnish a definition of the magnitude
of an infinite space, several methods for passing from finite to infinite
immediately suggest themselves.  Meckes~\cite{MecPDM} has shown that they are
largely equivalent.  Using the most elementary such method, coupled with some
Fourier analysis, we produce evidence for the following conjectural principle:
\slogan{magnitude subsumes all the most important invariants of integral
geometry} 
(Section~\ref{sec:inf}).  The most basic instance of this principle is the
fact that a line segment of length $t$ has magnitude $1 + t/2$, enabling one
to recover length from magnitude.  Less basic is the notion of the
\emph{magnitude dimension} of a space $A$, defined as the growth of the
function $t \goesto \mg{tA}$; here $tA$ is $A$ scaled up by a factor of $t$,
and $\mg{tA}$ is its magnitude.  We show, for example, that a subset of
$\reals^N$ with positive measure has magnitude dimension $N$.  At the cutting
edge is the conjecture (first stated in~\cite{AMSES}) that for any convex
subset $A$ of Euclidean space, all of the intrinsic volumes of $A$ can be
recovered from the function $t \goesto \mg{tA}$.

Review sections provide the necessary background on both enriched categories
and integral geometry.  No expertise in category theory or integral geometry
is needed to read this paper.

\paragraph*{Related work}  The basic ideas of this paper were first written up
in a 2008 internet posting~\cite{NCMS}.  Several papers have already built on
this.  Leinster and Willerton~\cite{AMSES} studied the large-scale asymptotics
of the magnitude of subsets of Euclidean space, and stated the conjecture just
mentioned.  That conjecture was partly motivated by numerical evidence and
heuristic arguments found by Willerton~\cite{WillHCC}, who also proved results
on the magnitude of Riemannian manifolds~\cite{WillMSS}.
Leinster~\cite{METAMB} established magnitude as maximum diversity.
Meckes~\cite{MecPDM}, \latin{inter alia}, proved the equivalence of several
definitions of the magnitude of compact metric spaces, and by using more
subtle analytical methods than are used here, extended some of the results of
Section~\ref{sec:inf} below.  The magnitude of spheres is especially well
understood~\cite{AMSES,WillMSS,MecPDM}.

In the literature on quantifying biodiversity, magnitude appears not only in
the paper of Solow and Polasky~\cite{SoPo}, but also in later papers such
as~\cite{POP}.  For an explanation of diversity in tune with the theory here,
see~\cite{MDISS}.

Geometry as the study of metric structures is developed in the books of
Blumenthal~\cite{BlumTAD} and Gromov~\cite{GroMSR}, among others;
representatives of the theory of finite metric spaces are~\cite{BlumTAD} and
papers of Dress and collaborators~\cite{BaDr,DHM}.  We will make contact with
the theory of spaces of negative type, which goes back to
Menger~\cite{MengMHR} and Schoenberg~\cite{Scho}.  This connection has been
exploited by Meckes~\cite{MecPDM}.  It is notable that the complete bipartite
graph $K_{3, 2}$ appears as a minimal example in both~\cite{BaDr} and
Example~\ref{eg:5-pt} below.

\paragraph*{Notation} Given $N \in \nat = \{0, 1, 2, \ldots\}$, we write
$\reals^N$ for real $N$-dimensional space as a set, topological space or
vector space---but with no implied choice of metric except when $N = 1$.  The
metric on a metric space $A$ is denoted by $d$ or $d_A$.  We write $\card X$
for the cardinality of a finite set $X$.  When $\cat{C}$ is a category, $C \in
\cat{C}$ means that $C$ is an object of $\cat{C}$.

\paragraph*{Acknowledgements}
I have had countless useful conversations on magnitude with Mark Meckes and
Simon Willerton.  Their insights have played a very important role, and I
thank them for it.  I also thank John Baez, Neal Bez, Paul Blackwell, Yemon
Choi, Christina Cobbold, David Corfield, Alastair Craw, Jacques Distler, Anton
Geraschenko, Martin Hyland, David Jordan, Andr\'e Joyal, Joachim Kock,
Christian Korff, Urs Schreiber, Josh Shadlen, Ivan Smith, David Speyer, and
Terry Tao.  Two web resources have
been crucial to the progress of this work: The $n$-Category Caf\'e%
\footnote{\href{http://golem.ph.utexas.edu/category}{http://golem.ph.utexas.edu/category}}
and MathOverflow.%
\footnote{\href{http://mathoverflow.net}{http://mathoverflow.net}}
Parts of this work were carried out
at the Centre de Recerca Matem\`atica (Barcelona) and the School of
Mathematics and Statistics at the University of Sheffield.  I thank them for
their hospitality.

\section{Enriched categories}
\label{sec:enr}

This section describes the conceptual origins of the notion of magnitude. 

We define the magnitude of an enriched category, in two steps.  First we
assign a number to every matrix; then we assign a matrix to every enriched
category.  We pause in between to recall some basic aspects of enriched
category theory: the definitions, and how a metric space can be viewed as an
enriched category.

\subsection{The magnitude of a matrix}
\label{subsec:mag-mx}


A \demph{rig} (or semiring) is a ri\emph{n}g without \emph{n}egatives: a set
$k$ equipped with a commutative monoid structure $(+, 0)$ and a monoid
structure $(\cdot, 1)$, the latter distributing over the former.  For us, rig
will mean \emph{commutative} rig: one whose multiplication is commutative.

It will be convenient to use matrices whose rows and columns are indexed by
abstract finite sets.  Thus, for finite sets $I$ and $J$, an $I \times J$
\demph{matrix} over a rig $k$ is a function $I \times J \to k$.  The usual
operations can be performed, e.g.\ an $H \times I$ matrix can be
multiplied by an $I \times J$ matrix to give an $H \times J$ matrix.  The
identity matrix is the Kronecker $\delta$.  An $I \times J$ matrix $\zeta$ has
a $J \times I$ transpose $\zeta^\tra$.  

Given a finite set $I$, we write $u_I \in k^I$ for the column vector with
$u_I(i) = 1$ for all $i \in I$.

\begin{defn}    \label{defn:mx-wtg}
Let $\zeta$ be an $I \times J$ matrix over a rig $k$.  A \demph{weighting} on
$\zeta$ is a column vector $w \in k^J$ such that $\zeta w = u_I$.  A
\demph{coweighting} on $\zeta$ is a row vector $v \in k^I$ such that $v \zeta
= u_J^\tra$.
\end{defn}

A matrix may admit zero, one, or many (co)weightings, but their freedom is
constrained by the following basic fact.

\begin{lemma}   \label{lemma:wtg-cowtg}
Let $\zeta$ be an $I \times J$ matrix over a rig, let $w$ be a weighting
on $\zeta$, and let $v$ be a coweighting on $\zeta$.  Then
\[
\sum_{j \in J} w(j) = \sum_{i \in I} v(i).
\]
\end{lemma}

\begin{proof}
$\sum_j w(j) 
= 
u_J^\tra w
=
v \zeta w
=
v u_I
=
\sum_i v(i)$.
\done
\end{proof}

We refer to the entries $w(j) \in k$ of a weighting $w$ as \demph{weights},
and similarly \demph{coweights}.  The lemma implies that if a matrix $\zeta$
has both a weighting and a coweighting, then the total weight is independent
of the weighting chosen.  This makes the following definition possible.

\begin{defn}    \label{defn:mx-mag}
A matrix $\zeta$ over a rig $k$ \demph{has magnitude} if it admits at least
one weighting and at least one coweighting.  Its \demph{magnitude} is then
\[
\mg{\zeta} = \sum_j w(j) = \sum_i v(i) \in k
\]
for any weighting $w$ and coweighting $v$ on $\zeta$.  
\end{defn}

We will be concerned with square matrices $\zeta$.  If
$\zeta$ is invertible then there are a unique weighting and a unique
coweighting.  (Conversely, if $k$ is a field then a unique weighting or
coweighting implies invertibility.)  The weights are then the sums of the rows
of $\zeta^{-1}$, and the coweights are the sums of the columns.
Lemma~\ref{lemma:wtg-cowtg} is obvious in this case, and there is an easy
formula for the magnitude:

\begin{lemma}   \label{lemma:invertible-basic}
Let $\zeta$ be an invertible $I \times I$ matrix over a rig.  Then $\zeta$
has a unique weighting $w$ given by $w(j) = \sum_i \zeta^{-1} (j, i)$ ($j \in
I$), and a unique coweighting given by the dual formula.  Also
\[
\mg{\zeta} = \sum_{i, j \in I} \zeta^{-1}(j, i).
\]
\ \done
\end{lemma}

Often our matrix $\zeta$ will be symmetric, in which case weightings and
coweightings are essentially the same.

\subsection{Background on enriched categories}
\label{subsec:enr-cats}

Here we review two standard notions: monoidal category, and category enriched
in a monoidal category.  

A \demph{monoidal category} is a category $\cat{V}$ equipped with an
associative binary operation $\otimes$ (which is formally a functor $\cat{V}
\times \cat{V} \to \cat{V}$) and a unit object $\Unit \in \cat{V}$.  The
associativity and unit axioms are only required to hold up to suitably
coherent isomorphism; see~\cite{MacCWM} for details.

\begin{examples}        \label{egs:mon-cats}
\begin{enumerate}
\item $\cat{V}$ is the category $\Set$ of sets, $\otimes$ is cartesian product
$\times$, and $\Unit$ is a one-element set $\{\star\}$.  
%

\item $\cat{V}$ is the category $\Vect$ of vector spaces over some field $K$,
the product $\otimes$ is the usual tensor product $\otimes_K$, and $\Unit =
K$. 

\item   \label{eg:mon-cats-reals}
A poset can be viewed as a category in which each hom-set has at most one
element.  In particular, consider the poset $([0, \infty], \geq)$ of
nonnegative reals together with infinity.  The objects of the resulting
category are the elements of $[0, \infty]$, there is one map $x \to
y$ when $x \geq y$, and there are none otherwise.  This is a monoidal category
with $\otimes = +$ and $\Unit = 0$.

\item   \label{eg:mon-cat-2}
Let $\Two$ be the category of Boolean truth values~\cite{LawMSG}: there are
two objects, $\fa$ (`false') and $\tr$ (`true'), and a single non-identity
map, $\fa \to \tr$.  Taking $\otimes$ to be conjunction and $\Unit = \tr$
makes $\Two$ monoidal.  Then $\Two$ is a monoidal subcategory
of $\Set$, identifying $\fa$ with $\emptyset$ and $\tr$ with $\{\star\}$.  It
is also a monoidal subcategory of $[0, \infty]$, identifying $\fa$ with
$\infty$ and $\tr$ with $0$.
\end{enumerate}
\end{examples}

Let $\cat{V} = (\cat{V}, \otimes, \Unit)$ be a monoidal category.  
The
definition of category enriched in $\cat{V}$, or $\cat{V}$-category, is
obtained from the definition of ordinary category by asking that the hom-sets
are no longer sets but objects of $\cat{V}$.  Thus, a (small)
\demph{$\cat{V}$-category} 
$\scat{A}$ consists of a set $\ob\scat{A}$ of
objects, an object $\Hom(a, b)$ of $\cat{V}$ for each $a, b \in \ob\scat{A}$,
and operations of composition and identity satisfying appropriate
axioms~\cite{KellBCE}.  The operation of composition consists of a map 
\[
\Hom(a, b) \otimes \Hom(b, c) \to \Hom(a, c)
\]
in $\cat{V}$ for each $a, b, c \in \ob\scat{A}$, while the identities are
provided by a map $\Unit \to \Hom(a, a)$ for each $a \in \ob\scat{A}$.  

There is an accompanying notion of enriched functor.  Given
$\cat{V}$-categories $\scat{A}$ and $\scat{A}'$, a \demph{$\cat{V}$-functor}
$F\cln \scat{A} \to \scat{A}'$ consists of a function $\ob\scat{A} \to
\ob\scat{A}'$, written $a \goesto F(a)$, together with a map
\[
\Hom(a, b) \to \Hom(F(a), F(b))
\]
in $\cat{V}$ for each $a, b \in \ob\scat{A}$, satisfying suitable
axioms~\cite{KellBCE}.  We write $\VCat$ for the category of 
$\cat{V}$-categories and $\cat{V}$-functors.

\begin{examples}        \label{egs:enr-cats}
\begin{enumerate}
\item Let $\cat{V} = \Set$.  Then $\VCat$ is the category $\Cat$ of (small)
categories and functors.

\item Let $\cat{V} = \Vect$.  Then $\VCat$ is the category of \demph{linear
categories} or \demph{algebroids}: categories equipped with a vector space
structure on each hom-set, such that composition is bilinear.

\item   \label{eg:enr-cats-reals}
Let $\cat{V} = [0, \infty]$.  Then, as observed by
Lawvere~\cite{LawMSG,LawTCS}, a $\cat{V}$-category is a \demph{generalized
metric space}.  That is, a $\cat{V}$-category consists of a set $A$ of objects
or points together with, for each $a, b \in A$, a real number $\Hom(a, b) =
d(a, b) \in [0, \infty]$, satisfying the axioms
\[
d(a, b) + d(b, c) \geq d(a, c),
\qquad
d(a, a) = 0
\]
($a, b, c \in A$).  Such spaces are more general than classical metric spaces
in three ways: $\infty$ is permitted as a distance, the separation axiom $d(a,
b) = 0 \implies a = b$ is dropped, and, most significantly, the symmetry axiom
$d(a, b) = d(b, a)$ is dropped.

A $\cat{V}$-functor $f\cln A \to A'$ between generalized metric spaces $A$ and
$A'$ is a \demph{distance-decreasing map}: one satisfying $ d(a, b) \geq
d(f(a), f(b)) $ for all $a, b \in A$.  Hence $[0, \infty]\hyph\Cat$ is the
category $\MS$ of generalized metric spaces and distance-decreasing maps.
Isomorphisms in $\MS$ are isometries.

\item Let $\cat{V} = \Two$.  A $\cat{V}$-category is a set equipped
with a preorder (a reflexive transitive relation), 
which up to equivalence of $\cat{V}$-categories is the same thing as a poset. 

The embedding $\Two \incl \Set$ of monoidal categories induces an embedding
$\Two\hyph\Cat \incl \Set\hyph\Cat$; this is the embedding $\Poset \incl \Cat$
of Example~\ref{egs:mon-cats}(\ref{eg:mon-cats-reals}).  Similarly,
the embedding $\Two \incl [0, \infty]$ induces an embedding $\Poset \incl
\MS$: as observed in~\cite{LawMSG}, a poset $(A, \leq)$ can be understood as a
non-symmetric metric space whose points are the elements of $A$ and whose
distances are all $0$ or $\infty$.
\end{enumerate}
\end{examples}

\subsection{The magnitude of an enriched category}
\label{subsec:mag-enr}

Here we meet the definition on which the rest of this work is built.

Having already defined the magnitude of a matrix, we now assign a
matrix to each enriched category.  To do this, we assume some further
structure on the enriching category $\cat{V}$.  In fact, we assume that
we have a notion of size for \emph{objects of} $\cat{V}$.  This, then, will
lead to a notion of size for \emph{categories enriched in} $\cat{V}$.

Let $\cat{V}$ be a monoidal category.  We will suppose given a rig $k$
and a monoid homomorphism
\[
\mg{\cdot}\cln
(\ob\cat{V}/\iso, \otimes, \Unit)
\to 
(k, \cdot, 1).
\]
(This is, deliberately, the same symbol as for magnitude; no confusion should
arise.)  The domain here is the monoid of isomorphism classes of objects of
$\cat{V}$.

\begin{examples}        \label{egs:size-invts}
\begin{enumerate}
\item   \label{eg:size-invt-Set} 
When $\cat{V}$ is the monoidal category $\Setf$ of finite sets, we take $k =
\rationals$ and $\mg{X} = \card{X}$.

\item When $\cat{V}$ is the monoidal category $\Vectf$ of finite-dimensional
vector spaces, we take $k = \rationals$ and $\mg{X} = \dim X$.

\item \label{eg:size-invt-R} When $\cat{V} = [0, \infty]$, we take $k =
\reals$ and $\mg{x} = e^{-x}$.  (If $\mg{\cdot}$ is to be measurable%
\footnote{I thank Mark Meckes for pointing out that the more obvious
hypothesis of continuity can be weakened to measurability~\cite{Frec}.  In
fact 
it is sufficient to assume that $\mg{\cdot}$
is bounded on some set of positive measure~\cite{Korm}.}
then the only possibility is $\mg{x} = C^x$ for some constant $C \geq 0$.)

\item When $\cat{V} = \Two$, we take $k = \integers$, $\mg{\fa} = 0$ and
$\mg{\tr} = 1$.  This is a restriction of the functions $\mg{\cdot}$
of~(\ref{eg:size-invt-Set}) and~(\ref{eg:size-invt-R}) along the
embeddings $\Two \incl \Setf$ and $\Two \incl [0, \infty]$ of
Example~\ref{egs:mon-cats}(\ref{eg:mon-cat-2}).
\end{enumerate}
\end{examples}

Write $\VCatf$ (with a small `c') for the category whose objects are the
$\cat{V}$-categories with \emph{finite} object-sets and whose maps are the
$\cat{V}$-functors between them.

\begin{defn}
Let $\scat{A} \in \VCatf$.  
\begin{enumerate}
\item The \demph{similarity matrix} of $\scat{A}$ is the
$\ob\scat{A} \times \ob\scat{A}$ matrix $\zeta_\scat{A}$ over $k$ defined by
$
\zeta_\scat{A}(a, b) = \mg{\Hom(a, b)}
$
($a, b \in \scat{A}$).  

\item A \demph{(co)weighting} on $\scat{A}$ is a (co)weighting on
$\zeta_\scat{A}$.

\item $\scat{A}$ \demph{has magnitude} if $\zeta_\scat{A}$ does; its
\demph{magnitude} is then $\mg{\scat{A}} = \mg{\zeta_\scat{A}}$.

\item $\scat{A}$ \demph{has M\"obius inversion} if $\zeta_\scat{A}$ is
invertible; its \demph{M\"obius matrix} is then $\mu_\scat{A} =
\zeta_{\scat{A}}^{-1}$.
\end{enumerate}
\end{defn}

Magnitude is, then, a partially-defined function $\mg{\cdot}\cln \VCatf
\partialmap k$.

\begin{examples}        \label{egs:magnitudes}
\begin{enumerate}
\item   \label{eg:magnitudes-Set}

When $\cat{V} = \Setf$, we obtain a notion of the magnitude $\mg{\scat{A}} \in
\rationals$ of a finite category $\scat{A}$ \cite{ECC,ECCSDS}.  This is
also called the \demph{Euler characteristic} of $\scat{A}$ and 
written $\chi(\scat{A})$.  There are theorems relating it to the Euler
characteristic of topological spaces, graphs, posets and orbifolds,
the cardinality of sets, and the order of groups.

Very many finite categories have M\"obius inversion (and in particular, Euler
characteristic).  The M\"obius matrix $\mu_\scat{A}$ is a generalization of
Rota's M\"obius function for posets~\cite{Rota}, which in turn generalizes the
classical M\"obius function on integers.  See~\cite{ECC} for explanation.

\item Similarly, taking $\cat{V} = \Vectf$ gives an invariant $\chi(\scat{A})
= \mg{\scat{A}} \in \rationals$ of linear categories $\scat{A}$ with finitely
many objects and finite-dimensional hom-spaces.

\item Taking $\cat{V} = [0, \infty]$ gives the notion of the magnitude $\mg{A}
\in \reals$ of a (generalized) finite metric space $A$.  This is the main
subject of this paper.

\item Taking $\cat{V} = \Two$ gives a notion of the magnitude $\mg{A} \in
\integers$ of a finite poset $A$.  Under the name of Euler characteristic,
this goes back to Rota~\cite{Rota}; see~\cite{StaEC1} for a modern account.
It is \emph{always} defined.  Indeed, every poset has M\"obius inversion, and
the M\"obius matrix is the M\"obius function of Rota mentioned
in~(\ref{eg:magnitudes-Set}).  

We have noted that a poset can be viewed as a category, or alternatively as a
non-symmetric metric space.  The notions of magnitude are compatible: the
magnitude of a poset is the same as that of the corresponding category or
generalized metric space.

\item Let $\cat{V}$ be a category of topological spaces in which every object
has a well-defined Euler characteristic (e.g.\ finite CW-complexes).
Taking $\mg{X}$ to be the Euler characteristic of a space $X$, we
obtain a notion of the magnitude or Euler characteristic of a topologically
enriched category.
\end{enumerate}
\end{examples}

The definition of the magnitude of a $\cat{V}$-category $\scat{A}$ is
independent of the composition and identities in $\scat{A}$, so could equally
well be made in the generality of $\cat{V}$-graphs.  (A
\demph{$\cat{V}$-graph} $\scat{G}$ is a set $\ob\scat{G}$ of objects together
with, for each $a, b \in \ob\scat{G}$, an object $\Hom(a, b)$ of $\cat{V}$.)
However, it is not clear that it is fruitful to do so.  
Two theorems on the magnitude or Euler characteristic of ordinary
categories, both proved in~\cite{ECC}, illuminate the general situation.

The first concerns directed graphs.  The Euler characteristic of a category
$\scat{A}$ is \emph{not} in general equal to the Euler characteristic of its
underlying graph $U(\scat{A})$.  But the functor $U$ has a left adjoint $F$,
assigning to a graph $\scat{G}$ the category $F(\scat{G})$ whose objects are
the vertices and whose maps are the paths in $\scat{G}$.  If $\scat{G}$ is
finite and circuit-free then $F(\scat{G})$ is finite, and the theorem is that
$\chi(F(\scat{G})) = \chi(\scat{G})$.  So the Euler characteristics of
categories and graphs are closely related, but not in the most obvious way.

The second theorem concerns the classifying space $B\scat{A}$ of a category
$\scat{A}$ (the geometric realization of its simplicial nerve).
Under suitable hypotheses, the topological space $B\scat{A}$ has a
well-defined Euler characteristic, and it is a theorem that $\chi(B\scat{A}) =
\chi(\scat{A})$.  It follows that if two categories have the same underlying
graph but different compositions then their classifying spaces, although not
usually homotopy equivalent, have the same Euler characteristic.  So if we
wish the Euler characteristic of a category to be defined in such a way that
it is equal to the Euler characteristic of its classifying space, it is
destined to be independent of composition.

\subsection{Properties}
\label{subsec:enr-props}

Much of ordinary category theory generalizes smoothly to enriched categories.
This includes many of the properties of the Euler characteristic of
categories~\cite{ECC}.  We list some of those properties now, using the
symbols $\cat{V}$, $k$ and $\mg{\cdot}$ as in the previous section.

There are notions of \demph{adjunction} and \demph{equivalence} between
$\cat{V}$-categories \cite{KellBCE}, generalizing the case $\cat{V} = \Set$ of
ordinary categories.  We write $\eqv$ for equivalence of $\cat{V}$-categories.

\begin{propn}   \label{propn:enr-eqv-adj}
Let $\scat{A}, \scat{B} \in \VCatf$.
\begin{enumerate}
\item   \label{part:enr-adj}
If there exist adjoint $\cat{V}$-functors $\scat{A} \oppairu \scat{B}$, and
$\scat{A}$ and $\scat{B}$ have magnitude, then $\mg{\scat{A}} =
\mg{\scat{B}}$.
\item   \label{part:enr-eqv-both}
If $\scat{A} \eqv \scat{B}$, and $\scat{A}$
and $\scat{B}$ have magnitude, then $\mg{\scat{A}} = \mg{\scat{B}}$.
\item   \label{part:enr-eqv-iff}
If $\scat{A} \eqv \scat{B}$ and $n \cdot 1 \in k$ has a multiplicative inverse
for all positive integers $n$, then $\scat{A}$ has magnitude if and only if
$\scat{B}$ does.
\end{enumerate}
\end{propn}

\begin{proof}
Part~(\ref{part:enr-adj}) has the same proof as Proposition~2.4(a)
of~\cite{ECC}, and part~(\ref{part:enr-eqv-both}) follows immediately.
Part~(\ref{part:enr-eqv-iff}) has the same proof as Lemma~1.12 of~\cite{ECC}.
\done
\end{proof}

For example, take a generalized metric space $A$ and adjoin a new point at
distance zero from some existing point.  Then the new space $A'$ is equivalent
to $A$.  By Proposition~\ref{propn:enr-eqv-adj}, if $A$ has magnitude then
$A'$ does too, and $\mg{A} = \mg{A'}$.  However, the proposition is trivial
for classical metric spaces $A, B$: if there is an adjunction between $A$ and
$B$ (and in particular if $A \eqv B$) then in fact $A$ and $B$ are isometric.

So far we have not used the multiplicativity of the function $\mg{\cdot}$ on
objects of $\cat{V}$.  We now show that it implies a multiplicativity
property of the function $\mg{\cdot}$ on $\cat{V}$-categories.

Assume that the monoidal category $\cat{V}$ is \demph{symmetric}, that is,
equipped with an isomorphism $X \otimes Y \to Y \otimes X$ for each pair
$X, Y$ of objects, satisfying axioms~\cite{MacCWM}.  There is a product on
$\VCat$, also denoted by $\otimes$, defined as follows.  Let $\scat{A},
\scat{B} \in \VCat$.  Then $\scat{A} \otimes \scat{B}$ is the
$\cat{V}$-category whose object-set is $\ob\scat{A} \times \ob\scat{B}$ and
whose hom-objects are given by
\[
\Hom((a, b), (a', b'))
=
\Hom(a, a') \otimes \Hom(b, b').
\]
Composition is defined with the aid of the symmetry~\cite{KellBCE}.  The unit
for this product is the one-object $\cat{V}$-category $\fcat{I}$ whose single
hom-object is $\Unit \in \cat{V}$.

\begin{examples}        \label{egs:enr-prods}
\begin{enumerate}
\item When $\cat{V} = \Set$, this is the ordinary product $\times$ of
categories.
\item   \label{eg:enr-prods-reals}
There is a family of products on metric spaces.  For $1 \leq p \leq
\infty$ and metric spaces $A$ and $B$, let $A \otimes_p B$ be the metric space
whose point-set is the product of the point-sets of $A$ and $B$, with
distances given by
\[
d((a, b), (a', b')) 
=
\begin{cases}
\bigl(d(a, a')^p + d(b, b')^p\bigr)^{1/p} &\text{if } p < \infty  \\
\max\{d(a, a'), d(b, b')\}      &\text{if } p = \infty.
\end{cases}
\]
Then the tensor product $\otimes$ defined above is $\otimes_1$.  
\end{enumerate}
\end{examples}

\begin{propn}   \label{propn:enr-prods}
Let $\scat{A}, \scat{B} \in \VCatf$.  If $\scat{A}$ and $\scat{B}$ have
magnitude then so does $\scat{A} \otimes \scat{B}$, with
\[
\mg{\scat{A} \otimes \scat{B}} = \mg{\scat{A}}\mg{\scat{B}}.
\]
Furthermore, the unit $\cat{V}$-category $\fcat{I}$ has magnitude $1$.
\end{propn}

\begin{proof}
As for Proposition~2.6 of~\cite{ECC}.
\done
\end{proof}

Magnitude is therefore a partially-defined monoid homomorphism
\[
\mg{\cdot}\cln 
(\VCatf/\eqv, \otimes, \fcat{I})
\partialmap
(k, \cdot, 1).
\]

Under mild assumptions, coproducts of $\cat{V}$-categories exist and interact
well with magnitude.  Indeed, assume that $\cat{V}$ has an initial object $0$,
with $X \otimes 0 \iso 0 \iso 0 \otimes X$ for all $X \in \cat{V}$.  Then for
any two $\cat{V}$-categories $\scat{A}$ and $\scat{B}$, the coproduct
$\scat{A} + \scat{B}$ in $\VCat$ exists.  It is constructed by taking the
disjoint union of $\scat{A}$ and $\scat{B}$ and setting $\Hom(a, b) = \Hom(b,
a) = 0$ whenever $a \in \scat{A}$ and $b \in \scat{B}$.  There is also an
initial $\cat{V}$-category $\emptyset$, with no objects.

When $\cat{V} = [0, \infty]$, the coproduct of metric spaces $A$ and $B$ is
their \demph{distant union}, the disjoint union of $A$ and $B$ with $d(a, b) =
d(b, a) = \infty$ whenever $a \in A$ and $b \in B$.

Assume also that $\mg{0} = 0$, where the $0$ on the left-hand side is the
initial object of $\cat{V}$.  This assumption and the previous ones hold in
all of our examples.

\begin{propn}
Let $\scat{A}, \scat{B} \in \VCatf$.  If $\scat{A}$ and $\scat{B}$ have
magnitude then so does $\scat{A} + \scat{B}$, with
\[
\mg{\scat{A} + \scat{B}} = \mg{\scat{A}} + \mg{\scat{B}}.
\]
Furthermore, the initial $\cat{V}$-category $\emptyset$ has magnitude $0$.
\end{propn}

\begin{proof}
As for Proposition~2.6 of~\cite{ECC}.
\done
\end{proof}

It might seem unsatisfactory that not every $\cat{V}$-category with finite
object-set has magnitude.  This can be resolved as follows.  

There are evident notions of algebra for a rig $k$ and (co)weighting for a
$\cat{V}$-category in a prescribed $k$-algebra.  As in
Lemma~\ref{lemma:wtg-cowtg}, the total weight is always equal to the total
coweight.  Given $\scat{A} \in \VCatf$, let $R(\scat{A})$ be the free
$k$-algebra containing a weighting $w$ and a coweighting $v$ for $\scat{A}$.
Then $\sum_a w(a) = \sum_a v(a) = \mgfml{\scat{A}}$, say.  This is
\emph{always} defined, and we may call $\mgfml{\scat{A}} \in R(\scat{A})$ the
\emph{formal magnitude} of $\scat{A}$.

A homomorphism $\phi$ from $R(\scat{A})$ to another $k$-algebra $S$ amounts to
a weighting and a coweighting for $\scat{A}$ in $S$, and
$\phi(\mgfml{\scat{A}}) \in S$ is independent of the homomorphism $\phi$
chosen.  In particular, $\scat{A}$ has magnitude in the original sense if and
only if there exists a $k$-algebra homomorphism $\phi\cln R(\scat{A}) \to k$;
in that case, $\mg{\scat{A}} = \phi(\mgfml{\scat{A}})$ for any such $\phi$.

This may lead to a more conceptually satisfactory theory, but at a price: the
magnitudes of different categories lie in different rigs, complicating results
such as those of the present section.  In any case, we say no more about this
approach.

\section{Finite metric spaces}
\label{sec:fin}

The definition of the magnitude of a finite metric space is a special case of
the definition for enriched categories.  Its most basic properties are special
cases of general results.  But metric spaces have many features not possessed
by enriched categories in general.  By exploiting them, we uncover a rich
theory.

A crucial feature of metric spaces is that they can be rescaled.  When handed
a space, we gain more information about it by considering the magnitudes
of its rescaled brothers and sisters than by taking it in isolation.  This
information is encapsulated in the so-called magnitude function of the
space.

For some spaces, the magnitude function exhibits wild behaviour:
singularities, negative magnitude, and so on.  But for geometrically orthodox
spaces such as subsets of Euclidean space, it turns out to be rather tame.
This is because they belong to the important class of `positive
definite' spaces.  Positive definiteness will play a central role when we come
to extend the definition of magnitude from finite to infinite spaces.  It is
explored thoroughly in the paper of Meckes~\cite{MecPDM}, who also describes
its relationship with the classical notion of negative type.

The term \demph{metric space} will be used in its standard sense, except that
$\infty$ is permitted as a distance.  Many of our theorems do hold for the
generalized metric spaces of
Example~\ref{egs:enr-cats}(\ref{eg:enr-cats-reals}), with the same proofs; but
to avoid cluttering the exposition, we leave it to the reader to discern which.

Throughout, we use matrices whose rows and columns are indexed by abstract
finite sets (as in Section~\ref{subsec:mag-mx}).  The identity matrix is
denoted by $\delta$.

\subsection{The magnitude of a finite metric space}
\label{subsec:mag-fin-ms}

We begin by restating the definitions from Section~\ref{sec:enr}, without
reference to enriched categories.  Let $A$ be a finite metric space.  Its
\demph{similarity matrix} $\zeta_A \in \reals^{A \times A}$ is defined by
$\zeta_A(a, b) = e^{-d(a, b)}$ ($a, b \in A$).  A \demph{weighting} on $A$ is
a function $w\cln A \to \reals$ such that $\sum_b \zeta_A(a, b) w(b) = 1$ for
all $a \in A$.  The space $A$ \demph{has magnitude} if it admits at least one
weighting; its \demph{magnitude} is then $\mg{A} = \sum_a w(a)$ for any
weighting $w$, and is independent of the weighting chosen.

A finite metric space $A$ \demph{has M\"obius inversion} if $\zeta_A$ is
invertible.  Its \demph{M\"obius matrix} is then $\mu_A = \zeta_A^{-1}$.  In
that case, there is a unique weighting $w$ given by $w(a) = \sum_b \mu(a,
b)$, and $\mg{A} = \sum_{a, b} \mu_A(a, b)$
(Lemma~\ref{lemma:invertible-basic}).  A generic real square matrix is
invertible; consequently, most finite metric spaces have M\"obius inversion.

Here are some elementary examples.

\begin{examples}        \label{egs:fin-ms-basic}
\begin{enumerate}
\item The empty space has magnitude $0$, and the one-point space has magnitude
$1$.  
\item   \label{eg:basic-two}
Let $A$ be the space consisting of two points distance $d$ apart.  Then 
\[
\zeta_{A}
=
\begin{pmatrix}
1       &e^{-d} \\
e^{-d}  &1      
\end{pmatrix}.
\]
This is invertible, so $A$ has M\"obius inversion and its magnitude is the sum
of all four entries of $\mu_A = \zeta_A^{-1}$:
\[
\mg{A} = 1 + \tanh(d/2)
\]
(Fig.~\ref{fig:two-points}).
\begin{figure}
\centering
\setlength{\unitlength}{1em}
\begin{picture}(26,10.2)(-1,0)
\put(1,1){\includegraphics[width=24.7\unitlength]{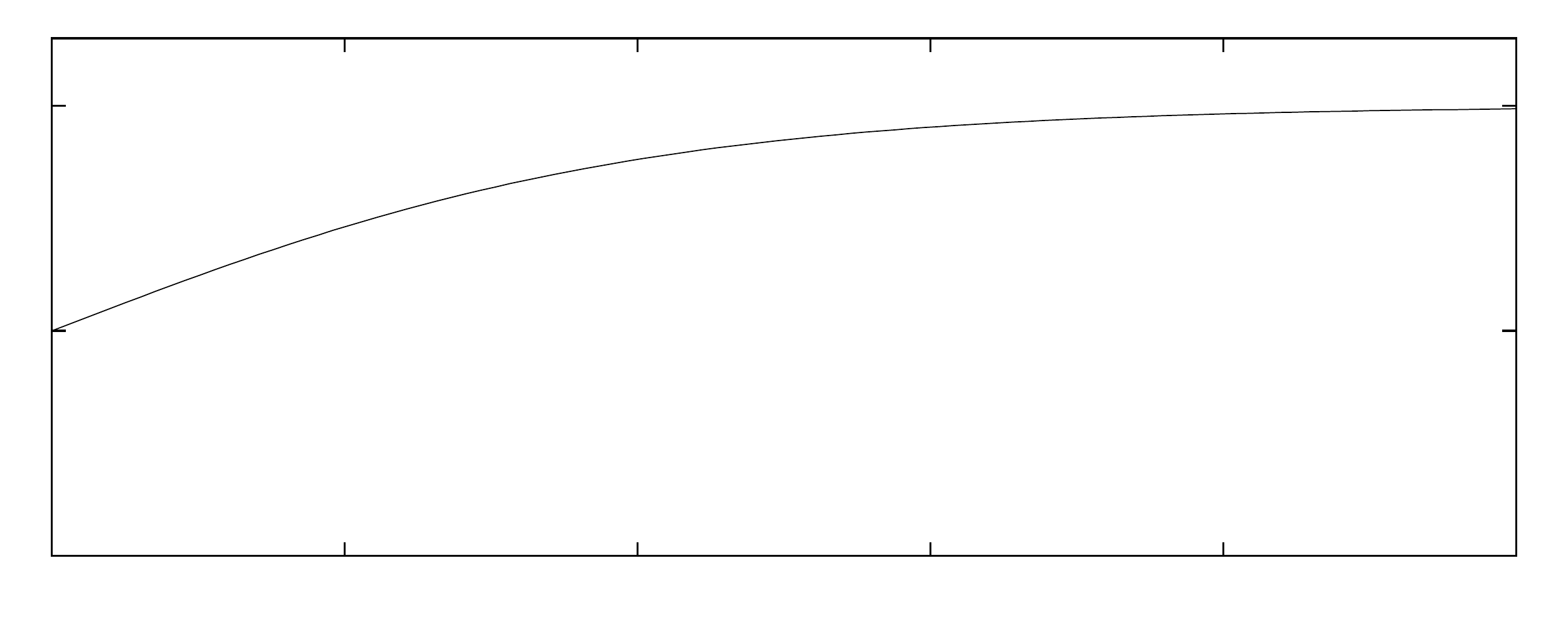}}
\qbezier[75](1.9,9.05)(13.4,9.05)(24.9,9.05)
\cell{1.9}{1}{b}{\scriptstyle{0}}
\cell{6.5}{1}{b}{\scriptstyle{1}}
\cell{11.1}{1}{b}{\scriptstyle{2}}
\cell{15.7}{1}{b}{\scriptstyle{3}}
\cell{20.3}{1}{b}{\scriptstyle{4}}
\cell{24.9}{1}{b}{\scriptstyle{5}}
\cell{1.5}{2}{r}{\scriptstyle{0}}
\cell{1.5}{5.5}{r}{\scriptstyle{1}}
\cell{1.5}{9}{r}{\scriptstyle{2}}
\cell{13.4}{0}{b}{d}
\cell{.5}{5.5}{r}{\mg{A}}
\end{picture}%
\caption{The magnitude of a two-point space}
\label{fig:two-points}
\end{figure}
This can be interpreted as follows.  When $d$ is small, $A$ closely resembles a
$1$-point space; correspondingly, the magnitude is little more than $1$.  As
$d$ grows, the points acquire increasingly separate identities and the
magnitude increases.  In the extreme, when $d = \infty$, the two points are
entirely separate and the magnitude is $2$.

\item   \label{eg:basic-discrete} 
A metric space $A$ is \demph{discrete}~\cite{LawTCS} if $d(a, b) = \infty$ for
all $a \neq b$ in $A$.  Let $A$ be a finite discrete space.  Then $\zeta_A$ is
the identity matrix $\delta$, each point has weight $1$, and $\mg{A} =
\card{A}$.  
\end{enumerate}
\end{examples}

The definition of the magnitude of a metric space first appeared in a paper of
Solow and Polasky~\cite{SoPo}, although with almost no mathematical
development.  They called it the `effective number of species', since the
points of their spaces represented biological species and the distances
represented inter-species differences (e.g.\ genetic).  We might say that the
magnitude of a metric space is the `effective number of points'.  Solow and
Polasky also considered the magnitude of correlation matrices, making
connections with the statistical concept of effective sample size.

Three-point spaces have magnitude; the formula follows from the proof of
Proposition~\ref{propn:3-pt}.  Meckes \cite[Theorem~\mrk{3.6}]{MecPDM} has
shown that four-point spaces have magnitude.  But spaces with five or more
points need not have magnitude (Example~\ref{eg:5-pt}).

We now describe two classes of space for which the magnitude exists and is
given by an explicit formula.

\begin{defn}
A finite metric space $A$ is \demph{scattered} if $d(a, b) > \log((\card{A})
- 1)$ for all distinct points $a$ and $b$.  (Vacuously, the empty space and
one-point space are scattered.) 
\end{defn}

\begin{propn}   \label{propn:scattered-mag}
A scattered space has magnitude.  Indeed, any scattered space $A$ has M\"obius
inversion, with M\"obius matrix given by the infinite sum
\[
\mu_A(a, b)     
=
\sum_{k = 0}^\infty
\sum_{\ a = a_0 \neq \cdots \neq a_k = b}
(-1)^k 
\zeta_A(a_0, a_1) \cdots \zeta_A(a_{k - 1}, a_k).
\]
\end{propn}

The inner sum is over all $a_0, \ldots, a_k \in A$ such that $a_0 = a$,
$a_k = b$, and $a_{j - 1} \neq a_j$ whenever $1 \leq j \leq k$.  That a
scattered space has magnitude was also proved in \cite[Theorem~2]{AMSES}, by a
different method that does not produce a formula for the M\"obius matrix.  

\begin{proof}
Write $n = \card{A}$.  For $a, b \in A$ and $k \geq 0$, put
\[
\mu_{A, k}(a, b)
=
\sum_{a = a_0 \neq \cdots \neq a_k = b}
\zeta_A(a_0, a_1) \cdots \zeta_A(a_{k - 1}, a_k).
\]
(In particular, $\mu_{A, 0}$ is the identity matrix.)  Write $\epsln =
\min_{a \neq b} d(a, b)$.  
Then
%
%
\begin{align*}        
\mu_{A, k + 1}(a, b)    &
=
\sum_{b'\cln b' \neq b}
\sum_{\ a = a_0 \neq \cdots \neq a_k = b'}
\zeta_A(a_0, a_1) \cdots \zeta_A(a_{k - 1}, b') \zeta_A(b', b)  \\
        &
\leq 
\sum_{b'\cln b' \neq b}
\sum_{\ a = a_0 \neq \cdots \neq a_k = b'}
\zeta_A(a_0, a_1) \cdots \zeta_A(a_{k - 1}, b') e^{-\epsln}     
=
e^{-\epsln}
\sum_{b' \cln b' \neq b}
\mu_{A, k}(a, b').
\end{align*}
The last sum is over $(n - 1)$ terms, so by induction, $\mu_{A, k}(a, b) \leq
\bigl( (n - 1) e^{-\epsln} \bigr)^k$ for all $a, b \in A$ and $k \geq 0$.
But $A$ is scattered, so $(n - 1) e^{-\epsln} < 1$, so the sum $\sum_{k =
0}^\infty (-1)^k \mu_{A, k}(a, b)$ converges for all $a, b \in A$.  A
telescoping sum argument finishes the proof.
%
\done
\end{proof}

\begin{defn}
A metric space is \demph{homogeneous} if its isometry group acts transitively
on points.
\end{defn}

\begin{propn}[Speyer~\cite{Spey}]        \label{propn:Speyer}
Every homogeneous finite metric space has magnitude.  Indeed, if
$A$ is a homogeneous space with $n \geq 1$ points then 
\[
\mg{A}
=
\frac{n^2}{\sum_{a, b} e^{-d(a, b)}}
=
\frac{n}{\sum_a e^{-d(x, a)}}
\]
for any $x \in A$.  There is a weighting $w$ on $A$ given by $w(a) = \mg{A}/n$
for all $a \in A$.
\end{propn}

\begin{proof}
By homogeneity, the sum $S = \sum_a \zeta_A(x, a)$ is independent of $x \in
A$.  Hence there is a weighting $w$ given by $w(a) = 1/S$ for all $a \in A$.
\done
\end{proof}

\begin{example} \label{eg:hgs-fin-comp}
For any (undirected) graph $G$ and $t \in (0, \infty]$, there is a metric
space $tG$ whose points are the vertices and whose distances are minimal
path-lengths, a single edge having length $t$.
%
%
Write $K_n$ for the complete graph on $n$ vertices.  Then
\[
\mg{tK_n} = \frac{n}{1 + (n - 1)e^{-t}}.
\]
\end{example}

In general, $e^{-d(a, b)}$ can be interpreted as the similarity or
closeness of the points $a, b \in A$ \cite{MDISS,SoPo}.
Proposition~\ref{propn:Speyer} states that the magnitude of a homogeneous
space is the reciprocal mean similarity.  

\begin{example} \label{eg:hgs-fin-bicomp}
A subspace can have greater magnitude than the whole space.  Let $K_{n, m}$ be
the graph with vertices $a_1, \ldots, a_n, b_1, \ldots, b_m$ and one edge
between $a_i$ and $b_j$ for each $i$ and $j$.  If $n$ is large then the mean
similarity between two points of $tK_{n, n}$ is approximately
$\frac{1}{2}(e^{-t} + e^{-2t})$ (Fig.~\ref{fig:K}).
\begin{figure}
\begin{center}
\setlength{\unitlength}{.8em}
\begin{picture}(18,11)
\put(0.2,0){\includegraphics[height=8.8em]{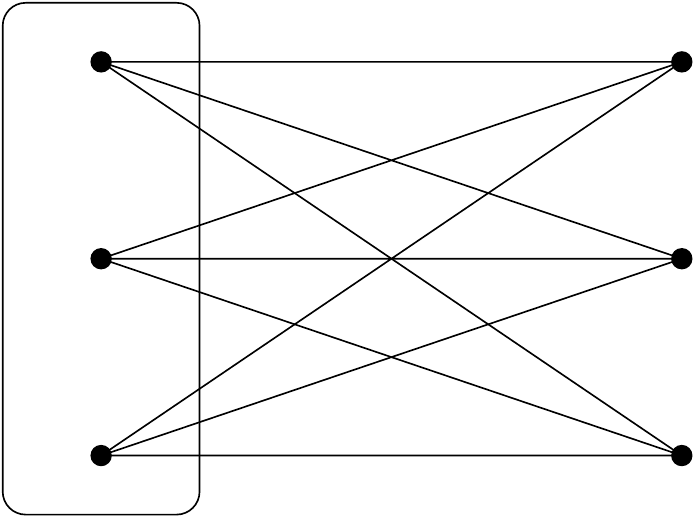}}
\cell{8.3}{0}{b}{t}
\cell{0.5}{9.6}{l}{a_1}
\cell{0.5}{5.4}{l}{a_2}
\cell{0.5}{1.2}{l}{a_3}
\cell{15.5}{9.6}{l}{b_1}
\cell{15.5}{5.4}{l}{b_2}
\cell{15.5}{1.2}{l}{b_3}
\end{picture}%
\end{center}
\caption{$tK_{n, n}$ and its subspace $2tK_n$, shown for $n = 3$.}
\label{fig:K} 
\end{figure}
%
On the other hand, $tK_{n, n}$ has a subspace $2tK_n = \{a_1,
\ldots, a_n\}$ in which the mean similarity is approximately $e^{-2t}$.  Since
$e^{-t} > e^{-2t}$, the mean similarity between points of $tK_{n, n}$ is
greater than that of its subspace $2tK_n$; hence $\mg{tK_{n, n}} <
\mg{2tK_n}$.  In fact, it can be shown using Proposition~\ref{propn:Speyer}
that $\mg{tK_{n, n}} < \mg{2tK_n}$ whenever $n > e^t + 1$.
\end{example}

\subsection{Magnitude functions}
\label{subsec:mag-fns}

In physical situations, distance depends on the choice of unit of length;
making a different choice rescales the metric by a constant factor.  In the
definition of $\mg{x}$ as $e^{-x}$
(Example~\ref{egs:size-invts}(\ref{eg:size-invt-R})), the constant $e^{-1}$
was chosen without justification; choosing a different constant between $0$
and $1$ also amounts to rescaling the metric.  For both these reasons, every
metric space should be seen as a member of the one-parameter family of spaces
obtained by rescaling it.

\begin{defn}
Let $A$ be a metric space and $t \in (0, \infty)$.  Then $tA$ denotes the
metric space with the same points as $A$ and
$
d_{tA}(a, b) = t d_A(a, b)
$
($a, b \in A$).
\end{defn}

Most familiar invariants of metric spaces behave in a predictable way when the
space is rescaled.  This is true, for example, of topological invariants,
diameter, and Hausdorff measure of any dimension.  But magnitude does not
behave predictably under rescaling.  Graphing $\mg{tA}$ against $t$ therefore
gives more information about $A$ than is given by $\mg{A}$ alone.

\begin{defn}
Let $A$ be a finite metric space.  The \demph{magnitude function} of $A$ is
the partially-defined function $t \mapsto \mg{tA}$, defined for all $t \in (0,
\infty)$ such that $tA$ has magnitude.
\end{defn}

\begin{examples}
\begin{enumerate}
\item Let $A$ be the space consisting of two points distance $d$ apart.  By
Example~\ref{egs:fin-ms-basic}(\ref{eg:basic-two}), the magnitude function of
$A$ is defined everywhere and given by $t \goesto 1 + \tanh(dt/2)$.

\item Let $A = \{a_1, \ldots, a_n\}$ be a nonempty homogeneous space, and
write $E_i = d(a_1, a_i)$.  By Proposition~\ref{propn:Speyer}, the magnitude
function of $A$ is
\[
t \goesto
n \Bigl/ \sum_{i = 1}^n e^{-E_i t}.
\]
In the terminology of statistical mechanics, the denominator is the partition
function for the energies $E_i$ at inverse temperature $t$.%
\footnote{I thank Simon Willerton for suggesting that some such relationship
should exist.}

\item Let $R$ be a finite commutative ring.  For $a \in R$, write
\[
\nu(a)
=
\min \{ k \in \nat \such a^{k + 1} = 0 \}
\in
\nat \cup \{\infty\}.
\]
There is a metric $d$ on $R$ given by $d(a, b) = \nu(b - a)$, and the
resulting metric space $A_R$ is homogeneous.  Write $q = e^{-t}$, and
$\Nil(R)$ for the set of nilpotent elements.  By
Proposition~\ref{propn:Speyer}, $A_R$ has magnitude function
%
\[
t 
\goesto
\mg{tA_R}
=
\card{R}
\Bigl/\!\!
\sum_{a \in \Nil(R)} q^{\nu(a)} 
=
\card{R}
\Bigl/
(1 - q) 
\sum_{k = 0}^\infty 
\card{ \{ a \in R \such a^{k + 1} = 0 \} } \cdot
q^k
\]
%
where the last expression is an element of the field $\rationals(\!(q)\!)$ of
formal Laurent series.  
\end{enumerate}
\end{examples}

To establish the basic properties of magnitude functions, we need some
auxiliary definitions and a lemma.  A vector $v \in \reals^I$ is
\demph{positive} if $v(i) > 0$ for all $i \in I$, and \demph{nonnegative} if
$v(i) \geq 0$ for all $i \in I$.  Recall the definition of distance-decreasing
map from Example~\ref{egs:enr-cats}(\ref{eg:enr-cats-reals}).

\begin{defn}
A metric space $A$ is an \demph{expansion} of a metric space $B$ if there
exists a distance-decreasing surjection $A \to B$. 
\end{defn}

\begin{lemma}   \label{lemma:expansion-nn}
Let $A$ and $B$ be finite metric spaces, each admitting a nonnegative
weighting.  If $A$ is an expansion of $B$ then $\mg{A} \geq \mg{B}$. 
\end{lemma}

\begin{proof}
Take a distance-decreasing surjection $f\cln A \to B$.  Choose a right
inverse function $g\cln B \to A$ (not necessarily distance-decreasing).  Then
$\zeta_B(f(a), b) \geq \zeta_A(a, g(b))$ for all $a \in A$ and $b \in B$.  Let
$w_A$ and $w_B$ be nonnegative weightings on $A$ and $B$ respectively.  Then
\[
\mg{A}
=
\sum_{a, b} w_A(a) \zeta_B(f(a), b) w_B(b)
\geq
\sum_{a, b} w_A(a) \zeta_A(a, g(b)) w_B(b)
=
\mg{B},
\]
as required.
\done
\end{proof}

\begin{propn}   \label{propn:mf-basic}
Let $A$ be a finite metric space.  Then:
\begin{enumerate}
\item   \label{part:mf-inv}
$tA$ has M\"obius inversion (hence magnitude) for all but finitely many $t >
0$. 
\item   \label{part:mf-ana}
The magnitude function of $A$ is analytic at all $t > 0$ such that
$tA$ has M\"obius inversion.
\item   \label{part:mf-pos}
For $t \gg 0$, there is a unique, positive, weighting on $tA$.
\item   \label{part:mf-inc}
For $t \gg 0$, the magnitude function of $A$ is increasing.
\item   \label{part:mf-lim}
$\mg{tA} \to \card{A}$ as $t \to \infty$.
\end{enumerate}
\end{propn}

\begin{proof}
We use the space $\reals^{A \times A}$ of real $A \times A$ matrices, and
its open subset $\GL{A}{\reals}$ of invertible matrices.  We also use the
notions of weighting on, and magnitude of, a matrix
(Section~\ref{subsec:mag-mx}).  For $\zeta \in
\GL{A}{\reals}$, the unique weighting $w_\zeta$ on $\zeta$ and the magnitude
of $\zeta$ 
are given by
\begin{equation}        \label{eq:wtg-GL}
w_\zeta(a) 
= 
\sum_{b \in A} \zeta^{-1}(a, b)
=
\sum_{b \in A} (\adj \zeta)(a, b)/\det \zeta,
\qquad
\mg{\zeta} 
=
\sum_{a \in A} w_{\zeta}(a)
\end{equation}
($a \in A$), where $\adj$ denotes the adjugate.

For~(\ref{part:mf-inv}), first note that $\zeta_{tA} \to \delta \in
\GL{A}{\reals}$ as $t \to \infty$; hence $\zeta_{tA}$ is invertible for $t \gg
0$.  The matrix $\zeta_{tA} = (e^{-td(a, b)})$ is defined for all $t \in
\complexes$, and $\det\zeta_{tA}$ is analytic in $t$.  But $\det\zeta_{tA}
\neq 0$ for real $t \gg 0$, so by analyticity, $\det\zeta_{tA}$ has only
finitely many zeros in $(0, \infty)$.

Part~(\ref{part:mf-ana}) follows from equations~\eqref{eq:wtg-GL}.

For~(\ref{part:mf-pos}), each of the functions $\zeta \goesto w_\zeta(a)$ ($a
\in A$) is continuous on $\GL{A}{\reals}$ by~\eqref{eq:wtg-GL}.  But
$w_\delta(a) = 1$ for all $a \in A$,
so there is a
neighbourhood $U$ of $\delta$ in $\GL{A}{\reals}$ such that $w_\zeta(a) > 0$
for all $\zeta \in U$ and $a \in A$.  Since $\zeta_{tA} \to \delta$ as $t \to
\infty$, we have $\zeta_{tA} \in U$ for all $t \gg 0$.

Part~(\ref{part:mf-inc}) follows from part~(\ref{part:mf-pos}) and
Lemma~\ref{lemma:expansion-nn}.  

For~(\ref{part:mf-lim}), $\lim_{t \to \infty} \mg{tA} = \mg{\lim_{t \to \infty}
\zeta_{tA}} = \mg{\delta} = \card{A}$.
%
%
\done
\end{proof}

Part~(\ref{part:mf-inv}) implies that magnitude functions have only finitely
many singularities.  Proposition~\ref{propn:scattered-pd-pw} will provide an
explicit lower bound for parts~(\ref{part:mf-pos}) and~(\ref{part:mf-inc}).
Part~(\ref{part:mf-lim}) also appeared as Theorem~3 of~\cite{AMSES}.
 
Many natural conjectures about magnitude are disproved by the following
example.  Later we will see that
subspaces of Euclidean space are less prone to surprising behaviour.%
\footnote{`Our approach to general metric spaces bears the undeniable imprint
of early exposure to Euclidean geometry.  We just love spaces sharing a common
feature with $\reals^n$.'  (Gromov~\cite{GroMSR}, page xvi.)}

\begin{example} \label{eg:5-pt}
Fig.~\ref{fig:K32} shows the magnitude function of the space $K_{3, 2}$
defined in Example~\ref{eg:hgs-fin-bicomp}.  It is given by
\begin{figure}
\centering
\setlength{\unitlength}{1em}
\begin{picture}(23,15.2)(-2,0)
\put(1,1){\includegraphics[width=20.5\unitlength]{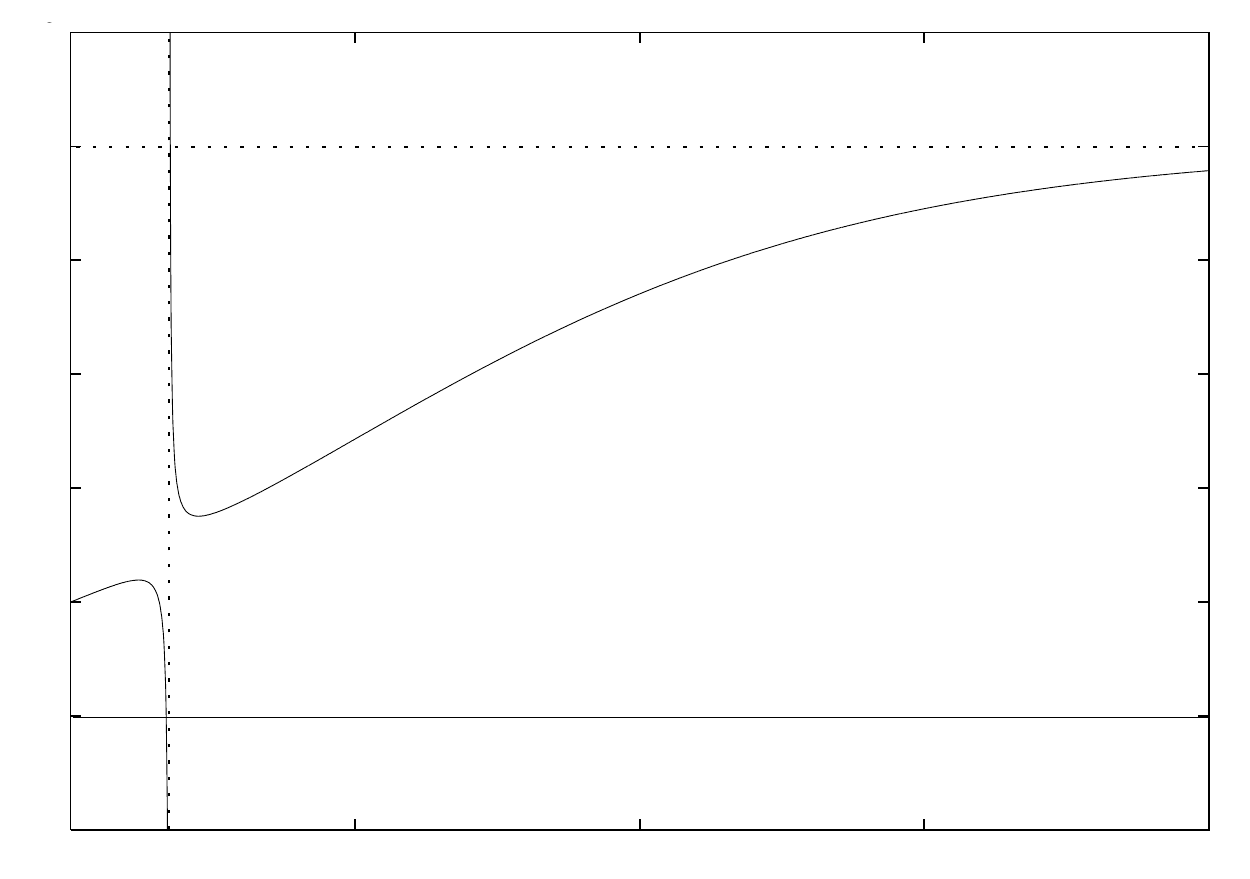}}
\cell{2.2}{1}{b}{\scriptstyle{0}}
\cell{6.85}{1}{b}{\scriptstyle{1}}
\cell{11.5}{1}{b}{\scriptstyle{2}}
\cell{16.15}{1}{b}{\scriptstyle{3}}
\cell{20.8}{1}{b}{\scriptstyle{4}}
\cell{4}{0.7}{b}{\scriptstyle{\log\sqrt{2}}}
\cell{1.8}{1.8}{r}{\scriptstyle{-1}}
\cell{1.8}{3.7}{r}{\scriptstyle{0}}
\cell{1.8}{5.5}{r}{\scriptstyle{1}}
\cell{1.8}{7.4}{r}{\scriptstyle{2}}
\cell{1.8}{9.3}{r}{\scriptstyle{3}}
\cell{1.8}{11.2}{r}{\scriptstyle{4}}
\cell{1.8}{13.0}{r}{\scriptstyle{5}}
\cell{1.8}{14.9}{r}{\scriptstyle{6}}
\cell{11.5}{-0.2}{b}{t}
\cell{-2}{8.4}{l}{\mg{tK_{3,2}}}
\end{picture}%
\caption{The magnitude function of the bipartite graph $K_{3, 2}$}
\label{fig:K32}
\end{figure}
\[
\mg{tK_{3,2}}
=
\frac{5 - 7e^{-t}}{(1 + e^{-t})(1 - 2e^{-2t})}
\]
($t \neq \log\sqrt{2}$); the magnitude of $(\log\sqrt{2})K_{3,2}$ is
undefined.  (One can compute this directly or use
Proposition~\ref{propn:const-dist}.)  Several features of the graph are
apparent.  At some scales, the magnitude is negative; at others, it is greater
than the number of points.  There are also intervals on which the magnitude
function is strictly decreasing.  Furthermore, this example shows that a space
with magnitude can have a subspace without magnitude: for
$(\log\sqrt{2})K_{3,2}$ is a subspace of $(\log\sqrt{2})K_{3,3}$, which, being
homogeneous, has magnitude (Proposition~\ref{propn:Speyer}).

(The graph $K_{3, 2}$ is also a well-known counterexample in the
theory of spaces of negative type~\cite{GrWi}.  The connection is
explained, in broad terms, by the remarks in Section~\ref{subsec:pd}.)
\end{example}

The first example of a finite metric space with undefined magnitude was found
by Tao~\cite{Tao}, and had $6$ points.  The first examples of $n$-point spaces
with magnitude outside the interval $[0, n]$ were found by the author and
Simon Willerton, and were again $6$-point spaces.

\begin{example} \label{eg:6-over-5}
This is an example of a space $A$ for which $\lim_{t \to 0} \mg{tA} \neq 1$,
due to Willerton (private communication, 2009).  Let $A$ be the graph $K_{3,
3}$ (Fig.~\ref{fig:K}) with three new edges adjoined: one from $b_i$ to $b_j$
whenever $1 \leq i < j \leq 3$.  Then $\mg{tA} = 6/(1 + 4e^{-t}) \to 6/5$ as
$t \to 0$.
\end{example}

\subsection{Constructions}

For each way of constructing a new metric space from old, we may ask whether
the magnitude of the new space is determined by the magnitudes of the old
ones.  Here we answer this question positively for four constructions: unions
(of a special type), tensor products, fibrations, and constant-distance
gluing.

\subsubsection*{Unions}

Let $X$ be a metric space with subspaces $A$ and $B$.  The magnitude of $A
\cup B$ is not in general determined by the magnitudes of $A$, $B$ and $A \cap
B$: consider one-point spaces.  In this respect, magnitude of metric spaces is
unlike cardinality of sets, for which there is the inclusion-exclusion
formula.  We do, however, have an inclusion-exclusion formula for magnitude
when the union is of a special type.

\begin{defn}    \label{defn:projects-to}
Let $X$ be a metric space and $A, B \sub X$.  Then $A$ \demph{projects to} $B$ 
if for all $a
\in A$ there exists $\pi(a) \in A \cap B$ such that for all $b \in B$,
\[
d(a, b) = d(a, \pi(a)) + d(\pi(a), b).
\]
\end{defn}

In this situation, $d(a, \pi(a)) = \inf_{b \in B} d(a, b)$.  If all distances
in $X$ are finite then $\pi(a)$ is unique for $a$.

\begin{propn}
\label{propn:union-fin}
Let $X$ be a finite metric space and $A, B \sub X$.  Suppose that $A$ projects
to $B$ and $B$ projects to $A$.  If $A$ and $B$ have magnitude then so does $A
\cup B$, with
\[
\mg{A \cup B} = \mg{A} + \mg{B} - \mg{A \cap B}.
\]
Indeed, if $w_A$, $w_B$ and $w_{A \cap B}$ are
weightings on $A$, $B$ and $A \cap B$ respectively then there is a weighting
$w$ on $A \cup B$ defined by 
\[
w(x)
=
\begin{cases}
w_A(x)  &\text{if } x \in A \without B  \\
w_B(x)  &\text{if } x \in B \without A  \\
w_A(x) + w_B(x) - w_{A \cap B}(x)       &\text{if } x \in A \cap B.
\end{cases}
\]
\end{propn}

\begin{proof}
Let $a \in A \without B$.  Choose a point $\pi(a)$ as in
Definition~\ref{defn:projects-to}.  Then
\begin{align*}
\sum_{x \in A \cup B} \zeta(a, x) w(x)  
        &
=
\sum_{a' \in A} \zeta(a, a') w_A(a')
+
\sum_{b \in B} \zeta(a, b) w_B(b)       
-
\sum_{c \in A \cap B} \zeta(a, c)w_{A \cap B}(c)        \\
        &
=
1
+
\zeta(a, \pi(a)) 
\biggl\{
\sum_{b \in B} \zeta(\pi(a), b) w_B(b)     
-
\sum_{c \in A \cap B} \zeta(\pi(a), c)w_{A \cap B}(c)
\biggr\}
=
1.
\end{align*}
Similar arguments apply when we start with a point of $B \without A$ or $A
\cap B$.  This proves that $w$ is a weighting, and the result follows.
\done
\end{proof}

It can similarly be shown that if $A$, $B$ and $A \cap B$ all have M\"obius
inversion then so does $A \cup B$.  The proof is left to the reader; we just
need the following special case.

\begin{cor}
\label{cor:one-point-fin}
Let $X$ be a finite metric space and $A, B \sub X$.  Suppose that $A \cap B$
is a singleton $\{c\}$, that for all $a \in A$ and $b \in B$,
\[
d(a, b) = d(a, c) + d(c, b),
\]
and that $A$ and $B$ have magnitude.  Then $A \cup B$ has magnitude $\mg{A} +
\mg{B} - 1$.  Moreover, if $A$ and $B$ have M\"obius inversion then so does $A
\cup B$, with
\[
\mu_{A \cup B}(x, y)
=
\begin{cases}
\mu_A(x, y)     &\text{if } x, y \in A \text{ and } (x, y) \neq (c, c) \\
\mu_B(x, y)     &\text{if } x, y \in B \text{ and } (x, y) \neq (c, c) \\
\mu_A(c, c) + \mu_B(c, c) - 1   &\text{if } (x, y) = (c, c) \\
0               &\text{otherwise.}
\end{cases}
\]
\end{cor}

\begin{proof}
The first statement follows from Proposition~\ref{propn:union-fin}, and
the second is easily checked.  
\done
\end{proof}

\begin{cor}     \label{cor:mag-reals}
Every finite subspace of $\reals$ has M\"obius inversion.  If $A = \{a_0 <
\cdots < a_n\} \sub \reals$ then, writing $d_i = a_i - a_{i - 1}$,
\[
\mg{A} 
= 
1 + \sum_{i = 1}^n \tanh\frac{d_i}{2}.
\]
The weighting $w$ on $A$ is given by
\[
w(a_i) = 
\frac{1}{2}
\left(
\tanh \frac{d_i}{2} + \tanh \frac{d_{i + 1}}{2}
\right)
\]
($0 \leq i \leq n$), where by convention $d_0 = d_{n + 1} = \infty$ and $\tanh
\infty = 1$.  
\end{cor}

\begin{proof}
This follows by induction from
Example~\ref{egs:fin-ms-basic}(\ref{eg:basic-two}),
Proposition~\ref{propn:union-fin} and Corollary~\ref{cor:one-point-fin}.  (An
alternative proof is given in \cite[Theorem~4]{AMSES}.)  
\done
\end{proof}

Thus, in a finite subspace of $\reals$, the weight of a point depends
only on the distances to its neighbours.  This is reminiscent of the Ising
model in statistical mechanics~\cite{DMS}, but whether there is any
substantial connection is unknown.

\begin{example}
The magnitude function is not a complete invariant of finite metric spaces.
Indeed, let $X = \{0, 1, 2, 3\} \sub \reals$.  Let $Y$ be the four-vertex
\textsf{Y}-shaped graph, viewed as a metric space as in
Example~\ref{eg:hgs-fin-comp}.  I claim that $X$ and $Y$ have the same
magnitude function, even though they are not isometric.  For put $A = \{0, 1,
2\} \sub \reals$ and $B = \{0, 1\} \sub \reals$.  Both $tX$ and $tY$ can be
expressed as unions, satisfying the hypotheses of
Corollary~\ref{cor:one-point-fin}, of isometric copies of $tA$ and $tB$.
Hence $ \mg{tX} = \mg{tA} + \mg{tB} - 1 = \mg{tY} $ for all $t > 0$.
\end{example}

\subsubsection*{Tensor products}

Recall from Example~\ref{egs:enr-prods}(\ref{eg:enr-prods-reals}) the
definition of the tensor product of metric spaces.
Proposition~\ref{propn:enr-prods} implies (and it is easy to prove directly):

\begin{propn}   \label{propn:prod-ms-fin}
If $A$ and $B$ are finite metric spaces with magnitude then $A
\otimes B$ has magnitude, given by $\mg{A \otimes B} = \mg{A} \mg{B}$.
\done
\end{propn}

\begin{example}
Let $q$ be a prime power, and denote by $\Field{q}$ the field of $q$ elements
metrized by $d(a, b) = 1$ whenever $a \neq b$.  Then for $N \in \nat$, the
metric tensor product $\Field{q}^{\otimes N}$ is the set $\Field{q}^N$ with
the Hamming metric.  Its magnitude function is
\[
t 
\goesto
\mg{t\Field{q}}^N
=
\left(
\frac{q}{1 + (q - 1)e^{-t}}
\right)^N
\]
by Example~\ref{eg:hgs-fin-comp} and Proposition~\ref{propn:prod-ms-fin}.

More generally, a \demph{linear code} is a vector subspace $C$ of
$\Field{q}^N$~\cite{MaSl}.  Its (single-variable) \demph{weight enumerator} is
the polynomial $W_C(x) = \sum_{i = 0}^N A_i(C) x^i \in \integers[x]$, where
$A_i(C)$ is the number of elements of $C$ whose Hamming distance from $0$ is
$i$.  Since $C$ is homogeneous, Proposition~\ref{propn:Speyer} implies that
its magnitude function is
\[
t \goesto (\card{C})/W_C(e^{-t}).
\]
The magnitude function of a code therefore carries the same, important,
information as its weight enumerator.
\end{example}

Similarly, if $A$ and $B$ are finite metric spaces with magnitude then their
coproduct or distant union $A + B$ (Section~\ref{subsec:enr-props}) has
magnitude $\mg{A + B} = \mg{A} + \mg{B}$.

\subsubsection*{Fibrations}

A fundamental property of the Euler characteristic of topological spaces is
its behaviour with respect to fibrations.  If a space $A$ is fibred over a
connected base $B$, with fibre $F$, then under suitable hypotheses, $\chi(A) =
\chi(B) \chi(F)$.  An analogous formula holds for the Euler characteristic of
a fibred category~\cite{ECC}.  

Apparently no general notion of fibration of enriched categories has yet been
formulated.  Nevertheless, we define here a notion of fibration of metric
spaces sharing common features with the categorical and topological notions,
and we prove an analogous theorem on magnitude.

\begin{figure}
\centering
\setlength{\unitlength}{1em}
\begin{picture}(13.5,10)(-1.5,0)
\put(0,5){\framebox(12,5)}
\cell{-.5}{7.5}{r}{A}
\put(0,1.5){\line(1,0){12}}
\cell{-.5}{1.5}{r}{B}
\cell{6}{3.25}{c}{\mbox{\Large\ensuremath{\downarrow}}}
\cell{6.5}{3.25}{l}{p}
\cell{3}{1.48}{c}{\zmark}  
\cell{3}{1.1}{t}{p(a)}
\cell{9}{1.48}{c}{\zmark}
\cell{9}{1.1}{t}{b'}
\cell{3}{7}{c}{\zmark}
\cell{2.5}{7}{r}{a}
\cell{9}{7}{c}{\zmark}
\cell{9.5}{7}{l}{a_{b'}}
\cell{9}{9}{c}{\zmark}
\cell{9.5}{9.25}{l}{a'}
\qbezier[20](9,5)(9,7.5)(9,10)
\end{picture}%
\caption{Fibration of metric spaces}
\label{fig:fib}
\end{figure}
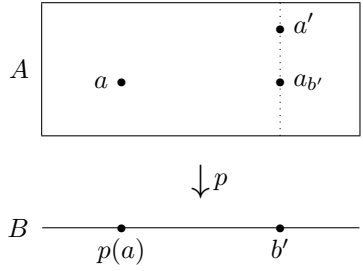

\begin{defn}
Let $A$ and $B$ be metric spaces.  A \demph{(metric) fibration} from $A$ to
$B$ is a distance-decreasing map $p\cln A \to B$ with the following property
(Fig.~\ref{fig:fib}): for all $a \in A$ and $b' \in B$ with $d(p(a), b') <
\infty$, there exists $a_{b'} \in p^{-1}(b')$ such that for all $a' \in
p^{-1}(b')$,
\begin{equation}        \label{eq:fib}
d(a, a')
=
d(p(a), b') + d(a_{b'}, a'). 
\end{equation}
\end{defn}

\begin{example}
Let $C_t$ be the circle of circumference $t$, metrized non-symmetrically by
taking $d(a, b)$ to be the length of the anticlockwise arc from $a$ to $b$.
(This is a generalized metric space in the sense of
Example~\ref{egs:enr-cats}(\ref{eg:enr-cats-reals}).)  Let $k$ be a positive
integer.  Then the $k$-fold covering $C_{kt} \to C_t$, locally an isometry, is
a fibration.
\end{example}

\begin{lemma}   \label{lemma:constant-fibre}
Let $p\cln A \to B$ be a fibration of metric spaces.  Let $b, b' \in B$ with
$d(b, b') < \infty$.  Then the fibres $p^{-1}(b)$ and $p^{-1}(b')$ are
isometric. 
\end{lemma}

\begin{proof}
Equation~\eqref{eq:fib} and finiteness of $d(b, b')$ imply that $a_{b'}$ is
unique for $a$, so we may define a function $\gamma_{b, b'}\cln p^{-1}(b) \to
p^{-1}(b')$ by $\gamma_{b, b'}(a) = a_{b'}$.  It is distance-decreasing: for
if $a, c \in p^{-1}(b)$ then
\[
d(b, b') + d(\gamma_{b, b'}(a), \gamma_{b, b'}(c)) 
= 
d(a, \gamma_{b, b'}(c)) 
\leq
d(a, c) + d(c, \gamma_{b, b'}(c))
=
d(a, c) + d(b, b'),
\]
giving $d(\gamma_{b, b'}(a), \gamma_{b, b'}(c)) \leq d(a, c)$ by finiteness of
$d(b, b')$.

There is a distance-decreasing map $\gamma_{b', b}\cln p^{-1}(b') \to
p^{-1}(b)$ defined in the same way.  It is readily shown that $\gamma_{b, b'}$
and $\gamma_{b', b}$ are mutually inverse; hence they are isometries.
\done
\end{proof}

Let $B$ be a nonempty metric space all of whose distances are finite, and let
$p\cln A \to B$ be a fibration.  The \demph{fibre} of $p$ is any of the
spaces $p^{-1}(b)$ ($b \in B$); it is well-defined up to isometry.

\begin{thm}     \label{thm:fib-mag}
Let $p\cln A \to B$ be a fibration of finite metric spaces.  Suppose that
$B$ is nonempty with $d(b, b') < \infty$ for all $b, b' \in B$, and that
$B$ and the fibre $F$ of $p$ both have magnitude.  Then $A$ has magnitude,
given by 
$
\mg{A} = \mg{B}\mg{F}.
$
\end{thm}

\begin{proof}
Choose a weighting $w_B$ on $B$.  Choose, for each $b \in B$, a weighting
$w_b$ on the space $p^{-1}(b)$.  For $a \in A$, put $w_A(a) =
w_{p(a)}(a)w_B(p(a))$.  It is straightforward to check that $w_A$ is a
weighting, and the theorem follows.
\done
\end{proof}

\begin{examples}
\begin{enumerate}
\item A trivial example of a fibration is a product-projection $B \otimes F
\to B$.  In that case, Theorem~\ref{thm:fib-mag} reduces to
Proposition~\ref{propn:prod-ms-fin}.

\item Let $B$ be a finite metric space in which the triangle inequality holds
strictly for every triple of distinct points.  Let $F$ be a finite metric
space of small diameter:
\[
\diam(F)
\leq
\min 
\bigl\{ 
d(b, b') + d(b', b'') - d(b, b'') 
\such 
b, b', b'' \in B,\ b \neq b' \neq b'' 
\bigr\}.
\]
Choose for each $b, b' \in B$ an isometry $\gamma_{b, b'}\cln F \to F$, in
such a way that $\gamma_{b, b}$ is the identity and $\gamma_{b', b} =
\gamma_{b, b'}^{-1}$.  Then the set $A = B \times F$ can
be metrized by putting
\[
d((b, c), (b', c')) 
=
d(b, b') + d(\gamma_{b, b'}(c), c')
\]
($b, b' \in B$, $c, c' \in F$).  The projection $A \to B$ is a fibration (but
not a product-projection unless $\gamma_{b', b''} \of \gamma_{b, b'} =
\gamma_{b, b''}$ for all $b, b', b''$).  So if $B$ and $F$ have magnitude,
$\mg{A} = \mg{B}\mg{F}$. 
\end{enumerate}
\end{examples}

Arguments similar to Lemma~\ref{lemma:constant-fibre} show that a fibration
over $B$ amounts to a family $(A_b)_{b \in B}$ of metric spaces together with
a distance-decreasing map $\gamma_{b, b'}\cln A_b \to A_{b'}$ for each $b, b'
\in B$ such that $d(b, b') < \infty$, satisfying the following three
conditions.  First, $\gamma_{b, b}$ is the identity for all $b \in B$.
Second, $\gamma_{b', b} = \gamma_{b, b'}^{-1}$.  Third,
\[
\sup_{a \in A_b}
d\bigl(
\gamma_{b', b''} \gamma_{b, b'} (a), \gamma_{b, b''}(a)
\bigr)
\leq
d(b, b') + d(b', b'') - d(b, b'')
\]
for all $b, b', b'' \in B$ such that $d(b, b'), d(b', b'') < \infty$.

\subsubsection*{Constant-distance gluing}

Given metric spaces $A$ and $B$ and a real number $D \geq \max\{\diam A, \diam
B\}/2$, there is a metric space $\codi{A}{D}{B}$ defined as follows.  As a
set, it is the disjoint union of $A$ and $B$.  The metric restricted to $A$ is
the original metric on $A$; similarly for $B$; and $d(a, b) = d(b, a) = D$ for
all $a \in A$ and $b \in B$.

\begin{propn}   \label{propn:const-dist}
Let $A$ and $B$ be finite metric spaces, and take $D$ as above.
Suppose that $A$ and $B$ have magnitude, with $\mg{A}\mg{B} \neq e^{2D}$.
Then $\codi{A}{D}{B}$ has magnitude
\[
\frac{\mg{A} + \mg{B} - 2e^{-D} \mg{A}\mg{B}}{1 - e^{-2D}\mg{A}\mg{B}}.
\]
\end{propn}

\begin{proof}
Given weightings $w_A$ on $A$ and $w_B$ on $B$, there is a weighting $w$ on
$\codi{A}{D}{B}$ defined by
\[
w(a)
=
\frac{1 - e^{-D}\mg{B}}{1 - e^{-2D}\mg{A}\mg{B}} w_A(a),
\qquad
w(b)
=
\frac{1 - e^{-D}\mg{A}}{1 - e^{-2D}\mg{A}\mg{B}} w_B(b)
\]
($a \in A$, $b \in B$).  The result follows.
\done
\end{proof}

This provides an easy way to compute the magnitude functions in
Examples~\ref{eg:5-pt} and~\ref{eg:6-over-5}.

\subsection{Positive definite spaces}
\label{subsec:pd}

We saw in Example~\ref{eg:5-pt} that the magnitude of a finite metric space
may be undefined, or smaller than the magnitude of one of its subspaces, or
even negative.  We now introduce a class of spaces for which no such behaviour
occurs.  Very many spaces of interest---including all subsets of Euclidean
space---belong to this class.

\begin{defn}    \label{defn:fin-pd}
A finite metric space $A$ is \demph{positive definite} if the matrix $\zeta_A$
is positive definite.
\end{defn}

We emphasize that positive definiteness of a matrix is meant in the
\emph{strict} sense.

\begin{lemma}   \label{lemma:pd-closure}
\begin{enumerate}
\item   \label{part:pd-has}
A positive definite space has M\"obius inversion.  
\item   \label{part:pd-tensor}
The tensor product of positive definite spaces is positive definite.
\item   \label{part:pd-sub}
A subspace of a positive definite space is positive definite.  
\end{enumerate}
\end{lemma}

\begin{proof}
Parts~(\ref{part:pd-has}) and~(\ref{part:pd-sub}) are elementary.
For~(\ref{part:pd-tensor}), $\zeta_{A \otimes B}$ is the Kronecker product
$\zeta_A \otimes \zeta_B$, and the Kronecker product of positive definite
matrices is positive definite. 
\done
\end{proof}

In particular, a positive definite space has magnitude and a unique weighting.

\begin{propn}   \label{propn:CS}
Let $A$ be a positive definite finite metric space.  Then
\[
\mg{A}
=
\sup_{v \neq 0} 
\frac{\bigl(\sum_{a \in A} v(a)\bigr)^2}%
{v^\tra \zeta_A v}
\]
where the supremum is over $v \in \reals^A \without\{0\}$ and $v^\tra$ denotes
the transpose of $v$.  A vector $v$ attains the supremum if and only if it is
a nonzero scalar multiple of the unique weighting on $A$.
\end{propn}


\begin{proof}
Since $\zeta_A$ is positive definite, we have the Cauchy--Schwarz inequality:
\[
(v^\tra \zeta_A v) \cdot (w^\tra \zeta_A w) 
\geq
(v^\tra \zeta_A w)^2
\]
for all $v, w \in \reals^A$, with equality if and only if one of $v$ and $w$
is a scalar multiple of the other.  Taking $w$ to be the unique weighting on
$A$ gives the result.
\done
\end{proof}

\begin{cor}     \label{cor:pd-sub-mag}
If $A$ is a positive definite finite metric space and $B \sub A$, then
$\mg{B} \leq \mg{A}$.
\done
\end{cor}

\begin{cor}
A nonempty positive definite finite metric space has magnitude $\geq 1$.
\done
\end{cor}

For any finite metric space $A$, the set $\Sing(A) = \{ t \in (0, \infty)
\such \zeta_{tA} \text{ is singular} \}$ is finite
(Proposition~\ref{propn:mf-basic}(\ref{part:mf-inv})).  When $\Sing(A) =
\emptyset$, put $\sup(\Sing(A)) = 0$.

\begin{propn}   \label{propn:eventually-pd}
Let $A$ be a finite metric space.  Then $tA$ is positive definite for all $t >
\sup(\Sing(A))$.  In particular, $tA$ is positive definite for all $t
\gg 0$.
\end{propn}

\begin{proof}
Write $\minev(\xi)$ for the minimum eigenvalue of a real symmetric $A \times
A$ matrix $\xi$.  Then $\minev(\xi)$ is continuous in $\xi$.  Also $\minev(\xi)
> 0$ if and only if $\xi$ is positive definite, and if $\minev(\xi) = 0$ then
$\xi$ is singular.

Now $\zeta_{tA} \to \delta$ as $t \to \infty$, and $\minev(\delta) = 1$, so
$\minev(\zeta_{tA}) > 0$ for all $t \gg 0$.  On the other hand,
$\minev(\zeta_{tA})$ is continuous and nonzero for $t > \sup(\Sing(A))$.
Hence $\minev(\zeta_{tA}) > 0$ for all $t > \sup(\Sing(A))$.
\done
\end{proof}

It follows that a space with M\"obius inversion at all scales also satisfies
an apparently stronger condition.

\begin{defn}
A finite metric space $A$ is \demph{stably positive definite} if $tA$ is
positive definite for all $t > 0$.
\end{defn}

\begin{cor}
Let $A$ be a finite metric space.  Then $tA$ has M\"obius inversion for all $t
> 0$ if and only if $A$ is stably positive definite.
\done
\end{cor}

\begin{example} \label{eg:mag-discts}
Let $A$ be the space of Example~\ref{eg:6-over-5}.  It is readily
shown that $tA$ has a unique weighting for all $t > 0$.  By the remarks after
Definition~\ref{defn:mx-mag}, $tA$ has M\"obius inversion
for all $t > 0$, so $A$ is stably positive definite.  Hence
magnitude is not continuous with respect to the Gromov--Hausdorff metric even
when restricted to stably positive definite finite spaces.
\end{example}

Meckes~\cite[Theorem~\mrk{3.3}]{MecPDM} has shown that a finite metric space
is stably positive definite if and only if it is of negative type.  By
definition, a finite metric space $A$ is of \demph{negative type} if $\sum_{a,
b} v(a) d(a, b) v(b) \leq 0$ for all $v \in \reals^A$ such that $\sum_a v(a) =
0$.  A general metric space $A$ is of negative type if every finite subspace
is of negative type, or equivalently if $(A, \sqrt{d_A})$ embeds isometrically
into some Hilbert space~\cite{Scho}.  Many important classes of space are
known to be of negative type, including those that we prove below to be stably
positive definite; see~\cite[Theorem~\mrk{3.6}]{MecPDM} for a list.  But
whereas the classical results on negative type tend to rely on embedding
theorems, we are able to bypass these and prove our results directly.

Lemma~\ref{lemma:expansion-nn} gave additional hypotheses on finite metric
spaces $A$ and $B$ guaranteeing that if $A$ is an expansion of $B$ then
$\mg{A} \geq \mg{B}$.  Some additional hypotheses are needed, since
not every magnitude function is increasing (Example~\ref{eg:5-pt}).  The
following will also do.

\begin{lemma}   \label{lemma:expansion-pd}
Let $A$ and $B$ be finite metric spaces.  Suppose that $A$ is positive
definite and $B$ admits a nonnegative weighting.  If $A$ is an expansion of
$B$ then $\mg{A} \geq \mg{B}$.
\end{lemma}

\begin{proof}
First consider a distance-decreasing \emph{bi}jection $f\cln A \to B$.  Choose
a nonnegative weighting $w_B$ on $B$.  Without loss of generality, $f$ is the
identity as a map of sets; thus, $\zeta_A(a, a') \leq \zeta_B(a, a')$ for all
points $a, a'$.  Hence
\[
\mg{A}
\geq
\frac{(\sum w_B(a))^2}{w_B^\tra \zeta_A w_B}
\geq 
\frac{(\sum w_B(a))^2}{w_B^\tra \zeta_B w_B}
=
\mg{B},
\]
by Proposition~\ref{propn:CS}.

Now consider the general case of a distance-decreasing surjection from $A$ to
$B$.  We may choose a subspace $A' \sub A$ and a distance-decreasing bijection
$A' \to B$.  The space $A'$ is positive definite, so $\mg{A'} \geq \mg{B}$ by
the previous argument; but also $\mg{A} \geq \mg{A'}$ by
Corollary~\ref{cor:pd-sub-mag}. 
\done
\end{proof}

A positive definite space cannot have negative magnitude, but the following
example shows that it can have magnitude greater than the number of points.

\begin{example}
Take the space $K_{3, 2}$ of Example~\ref{eg:5-pt}.  It is easily shown that
$\Sing(K_{3, 2}) = \{ \log\sqrt{2} \}$.  Choose $u > \log\sqrt{2}$ such that
$\mg{u K_{3, 2}} > 5$ (say, $u = 0.35$): then $A = uK_{3, 2}$ is positive
definite by Proposition~\ref{propn:eventually-pd}, and $\mg{A} > \card{A}$.

This example also shows that a positive definite expansion of a positive
definite space may have smaller magnitude: for if $s > 1$ then $sA$ is an
expansion of $A$, but $\mg{sA} < \mg{A}$ (Fig.~\ref{fig:K32}).
\end{example}

A different positivity condition is sometimes useful: the existence of a
nonnegative weighting.  

\begin{lemma}
Let $A$ be a finite metric space admitting a nonnegative weighting.  Then $0
\leq \mg{A} \leq \card{A}$.
\end{lemma}

\begin{proof}
Choose a nonnegative weighting $w$ on $A$.  For all $a \in A$ we have $0 \leq
w(a) \leq (\zeta_A w)(a) = 1$, so $0 \leq w(a) \leq 1$.
%
%
Summing, $0 \leq \mg{A} \leq \card{A}$.
\done
\end{proof}

We now list some classes of space that are positive definite, or have positive
weightings, or both.  The Euclidean case will be covered in the next section.

\begin{propn}
\label{propn:R-pd}
Every finite subspace of $\reals$ is positive definite with positive weighting.
\end{propn}

\begin{proof}
Let us temporarily say that a finite metric space $A$ is \emph{good} if it has
M\"obius inversion and for all $v \in \reals^A$,
\[
v^\tra \mu_A v
\geq
\max_{a \in A} v(a)^2.
\]
I claim that if $A \cup B$ is a union of the type in
Corollary~\ref{cor:one-point-fin} and $A$ and $B$ are both good, then $A \cup
B$ is good.  Indeed, let $v \in \reals^{A \cup B}$.  By
Corollary~\ref{cor:one-point-fin}, 
\[
v^\tra \mu_{A \cup B} v
=
v\restr{A}^\tra \mu_A v\restr{A}
+
v\restr{B}^\tra \mu_B v\restr{B}
-
v(c)^2 
\]
where $v\restr{A}$ is the restriction of $v$ to $A$.  Now let $x \in A \cup
B$.  Without loss of generality, $x \in A$.  Since $A$ is good,
$v\restr{A}^\tra \mu_A v\restr{A} \geq v(x)^2$.  Since $B$ is good,
$v\restr{B}^\tra \mu_B v\restr{B} \geq v(c)^2$.  Hence $v^\tra \mu_{A \cup B}
v \geq v(x)^2$, proving the claim.

Every metric space with $0$, $1$ or $2$ points is good.  Every finite subset
of $\reals$ with $3$ or more points can be expressed nontrivially as a union
of the type in Corollary~\ref{cor:one-point-fin}.  It follows by induction
that every finite subset of $\reals$ is good and therefore positive definite.

Positivity of the weighting is immediate from Corollary~\ref{cor:mag-reals}.
\done
\end{proof}

For $N \in \nat$ and $1 \leq p \leq \infty$, write $\ell_p^N =
\reals^{\otimes_p N}$, where $\otimes_p$ is as defined in
Example~\ref{egs:enr-prods}(\ref{eg:enr-prods-reals}).  Thus, $\ell_p^N$ is
$\reals^N$ with the metric induced by the $p$-norm, $\nrm{x}{p} = (\sum_r
|x_r|^p)^{1/p}$.   

\begin{thm}   \label{thm:l1-pd}
Every finite subspace of $\ell_1^N$ is positive definite.
\end{thm}

\begin{proof}
Let $A$ be a finite subspace of $\ell_1^N$.  Write $\pr_1, \ldots, \pr_N\cln
\ell_1^N \to \reals$ for the projections.  Each space $\pr_r A$ is positive
definite by Proposition~\ref{propn:R-pd}, so $\prod_{r = 1}^N \pr_r A \sub
\ell_1^N$ is positive definite by
Lemma~\ref{lemma:pd-closure}(\ref{part:pd-tensor}), so $A$ is positive
definite by Lemma~\ref{lemma:pd-closure}(\ref{part:pd-sub}).  
\done
\end{proof}

In the category $\MS$ of metric spaces and distance-decreasing maps
(Example~\ref{egs:enr-cats}(\ref{eg:enr-cats-reals})), the categorical product
$\times$ is $\otimes_\infty$.  The class of positive definite spaces is
\emph{not} closed under $\times$.  For if it were then, by an argument similar
to the proof of Theorem~\ref{thm:l1-pd}, every finite subspace of
$\ell_\infty^N$ would be positive definite.  But in fact, \emph{every} finite
metric space embeds isometrically into $\ell_\infty^N$ for some
$N$~\cite{Scho}, whereas not every finite metric space is positive definite.
Comprehensive results on (non-)preservation of positive definiteness by the
products $\otimes_p$ can be found in~\cite[Section~\mrk{3.2}]{MecPDM}.

\begin{propn}   \label{propn:3-pt}
Every space with $3$ or fewer points is positive definite with positive
weighting.   
\end{propn}

\begin{proof}
The proposition is trivial for spaces with $2$ or fewer points.  Now take a
$3$-point space $A = \{a_1, a_2, a_3\}$, writing $Z_{ij} = \zeta(a_i, a_j)$.
We use Sylvester's criterion: a symmetric real $n \times n$ matrix is positive
definite if and only if the upper-left $m \times m$ submatrix has positive
determinant whenever $1 \leq m \leq n$.  This holds for $Z$ when $m = 1$ or $m
= 2$, and 
\begin{align*}
\det Z  &
=      
(1 - Z_{12})(1 - Z_{23})(1 - Z_{31})    \\
        &       
\quad
+
(1 - Z_{12})(Z_{12} - Z_{13} Z_{32}) 
+
(1 - Z_{23})(Z_{23} - Z_{21} Z_{13})
+
(1 - Z_{31})(Z_{31} - Z_{32} Z_{21})    
\end{align*}
which is positive by the triangle inequality.
The unique weighting is $v/\det Z$, where
%
\[
v_1     
=      
(1 - Z_{12})(1 - Z_{23})(1 - Z_{31})
+
(1 - Z_{23})(Z_{23} - Z_{21}Z_{13})     
>      
0
\]
%
and similarly $v_2$ and $v_3$.
\done
\end{proof}

Meckes~\cite[Theorem~\mrk{3.6}]{MecPDM} has shown that $4$-point spaces are
also positive definite.  By Example~\ref{eg:5-pt}, his result is optimal.

\begin{example}
The weighting on a $4$-point space may have negative components, as may the
weighting on a finite subspace of $\ell_1^N$.  Indeed, using
Proposition~\ref{propn:union-fin} one can show that in the space $\{(0, 0),
(t, 0), (0, t), (-t, 0)\} \sub \ell_1^2$, the weight at $(0, 0)$ is negative
whenever $t < \log 2$.
%
\end{example}

\begin{propn}   \label{propn:scattered-pd-pw}
Every scattered space is positive definite with positive weighting.
\end{propn}

This provides an alternative, quantitative proof that every finite metric
space, when scaled up sufficiently, becomes positive definite with
positive weighting (Propositions~\ref{propn:mf-basic}
and~\ref{propn:eventually-pd}).

\begin{proof}
Let $A$ be a scattered space with $n \geq 2$ points.  For positive
definiteness, we use the same argument as appears in the proof
of~\cite[Theorem 2]{AMSES}.  Let $v \in \reals^A$.  Then
\begin{align*}
v^\tra \zeta_A v  &
= 
\sum_a v(a)^2 + \sum_{a \neq b} v(a) \zeta_A(a, b) v(b) 
\geq
\sum_a v(a)^2 - \frac{1}{n - 1} \sum_{a \neq b} |v(a)||v(b)|    \\
&
=
\frac{1}{2(n - 1)} \sum_{a \neq b} 
\bigl( |v(a)| - |v(b)| \bigr)^2 
\geq 
0.
\end{align*}
The inequality $\zeta_A(a, b) < 1/(n - 1)$ ($a \neq b$) is strict, so if
$v^\tra \zeta_A v = 0$ then $v = 0$.  

To show that the unique weighting $w_A$ on $A$ is positive, we use the
proof of Proposition~\ref{propn:scattered-mag}.  There we showed that $A$ has
M\"obius inversion and that the M\"obius matrix is a sum $\mu_A = \sum_{k =
0}^\infty (-1)^k \mu_{A, k}$, where the matrices $\mu_{A, k}$ satisfy
\begin{equation}        \label{eq:mu-ineq}
\mu_{A, k + 1}(a, b) 
< 
\frac{1}{n - 1}
\sum_{b'\cln b' \neq b} \mu_{A, k}(a, b')
\end{equation}
for all $a, b$.  Hence $w_A = \sum_{k = 0}^\infty (-1)^k w_{A, k}$, where
$w_{A, k}(a) = \sum_b \mu_{A, k}(a, b)$.  Summing~\eqref{eq:mu-ineq} over all
$b \in A$ gives
\[
w_{A, k + 1}(a) 
< 
\frac{1}{n - 1} \sum_{b, b'\cln b' \neq b} \mu_{A, k}(a, b')
=
w_{A, k}(a)
\]
($a \in A$).  Hence $w_A(a) = \sum_{k = 0}^\infty (-1)^k w_{A, k}(a) > 0$ for
all $a \in A$. 
\done
\end{proof}

A metric space $A$ is \demph{ultrametric} if
$
\max\{d(a, b), d(b, c)\} \geq d(a, c)
$
for all $a, b, c \in A$.  

\begin{propn}   \label{propn:um-pd}
Every finite ultrametric space is positive definite with positive weighting.
\end{propn}

Positive definiteness was proved by Varga and Nabben~\cite{VaNa}, and
positivity of the weighting (rather indirectly) by Pavoine, Ollier and
Pontier~\cite{POP}.  Another proof of positive definiteness is given by
Meckes~\cite[Theorem~\mrk{3.6}]{MecPDM}.  Both parts of the following proof
are different from those cited.

\begin{proof}
Let $\WUM$ be the set of symmetric matrices $Z$ over $[0, \infty)$ such that
$Z_{ik} \geq \min\{Z_{ij}, Z_{jk}\}$ for all $i, j, k$ and $Z_{ii} > \max_{j
\neq k} Z_{jk}$ for all $i$.  (For a $1 \times 1$ matrix, this
maximum is to be interpreted as $0$.)  We show by
induction that every matrix in $\WUM$ is positive definite and that its unique
weighting (Definition~\ref{defn:mx-wtg}) is positive.  The proposition will
follow immediately.

The result is trivial for $0 \times 0$ and $1 \times 1$ matrices.  Now let $Z
\in \WUM$ be an $n \times n$ matrix with $n \geq 2$.  Put
$z = \min_{i, j} Z_{ij}$.  There is an equivalence relation $\sim$ on $\{1,
\ldots, n\}$ defined by $i \sim j$ if and only if $Z_{ij} > z$.

It is not the case that $i \sim j$ for all $i, j$.  Hence we may partition
$\{1, \ldots, n\}$ into two nonempty subsets that are each a union of
equivalence classes: say $\{1, \ldots, m\}$ and $\{m + 1, \ldots n\}$.  We
have $Z_{ij} = z$ whenever $i \leq m < j$, so $Z$ is a block sum
\[
Z 
=
\begin{pmatrix}
Z'              &z\UM{m}{n-m}   \\
z\UM{n-m}{m}    &Z''
\end{pmatrix}
\]
where $\UM{k}{\ell}$ denotes the $k \times \ell$ matrix all of whose entries
are $1$.  Since $Z' \in \WUM$ and $Z'_{ij} = Z_{ij} \geq z$ for all $i, j \leq
m$, we have $Y' = Z' - z\UM{m}{m} \in \WUM$.  Similarly, $Y'' = Z'' - z\UM{n -
m}{n - m} \in \WUM$, and
\[
Z 
=
z\UM{n}{n}
+
\begin{pmatrix}
Y'      &0      \\
0       &Y''
\end{pmatrix}.
\]
The first summand is positive semidefinite.  By inductive hypothesis, $Y'$ and
$Y''$ are positive definite, so the second summand is positive definite.
Hence $Z$ is positive definite.

Also by inductive hypothesis, $Y'$ and $Y''$ have positive weightings $v'$ and
$v''$ respectively.  Let $v$ be the concatenation of $v'$ and $v''$.  It is
straightforward to verify that
\[
\frac{v}{z(\mg{Y'} + \mg{Y''}) + 1}
%
\]
is a weighting on $Z$, and it is positive since $v'$ and $v''$ are positive and
$z, \mg{Y'}, \mg{Y''} \geq 0$.
\done
\end{proof}

\begin{cor}     \label{cor:um-diam}
If $A$ is a finite ultrametric space then $\mg{A} \leq e^{\diam A}$.
\end{cor}

\begin{proof}
Let $\Delta$ be the metric space with the same point-set as $A$ and $d(a, b) =
\diam A$ for all distinct points $a, b$.  By Proposition~\ref{propn:Speyer},
$\mg{\Delta} \leq e^{\diam A}$ and $\Delta$ has a positive weighting.  But
$\Delta$ is an expansion of $A$, so $\mg{A} \leq \mg{\Delta}$ by
Lemma~\ref{lemma:expansion-nn}.  \done
\end{proof}

A homogeneous space always has a positive weighting, by
Proposition~\ref{propn:Speyer}.  However, Example~\ref{eg:hgs-fin-bicomp} and
Corollary~\ref{cor:pd-sub-mag} together show that a homogeneous space need not
be positive definite.  A homogeneous space need not even have M\"obius
inversion: $(\log 2)K_{3, 3}$ is an example.  In particular, a finite metric
space may have magnitude but not M\"obius inversion.

Magnitude can be understood in terms of entropy or diversity.  For every
finite metric space $A$ and $q \in [0, \infty]$, there is a function $D_q^A$
assigning to each probability distribution $\pv{p}$ on $A$ a real number
$D_q^A(\pv{p})$, the \demph{diversity of order $q$} of the
distribution~\cite{METAMB, MDISS}.  An ecological community can be modelled as
a metric space $A$ (as in Section~\ref{subsec:mag-fin-ms}) together with a
probability distribution $\pv{p}$ on $A$ (representing the relative abundances
of the species).  Then $D_q^A(\pv{p})$ is a measure of the biodiversity of the
community.  In the special case that $A$ is discrete, the diversities are the
exponentials of the R\'enyi entropies~\cite{Ren}, and in particular, the
diversity of order 1 is the exponential of Shannon entropy.

It is a theorem~\cite{METAMB} that for each finite metric space $A$, there is
some probability distribution $\pv{p}$ maximizing $D_q^A(\pv{p})$ for all $q
\in [0, \infty]$ simultaneously.  Moreover, the maximal value of
$D_q^A(\pv{p})$ is independent of $q$; call it $\Dmax{A}$.  If $A$ is positive
definite with nonnegative weighting then, in fact, $\mg{A} = \Dmax{A}$:
magnitude is maximum diversity.

\subsection{Subsets of Euclidean space}
\label{subsec:fin-euc}

Here we show that every finite subspace of Euclidean space $\ell_2^N$ is
positive definite.  In particular, every such space has well-defined
magnitude.

Write $L_1(\reals^N)$ for the space of Lebesgue-integrable complex-valued
functions on $\reals^N$.  Define the Fourier transform $\ft{f}$ of $f \in
L_1(\reals^N)$ by
\[
\ft{f}(\xi)
=
\int_{\reals^N} e^{-2\pi i \ip{\xi}{x}} f(x) \dee x
\]
($\xi \in \reals^N$).  Define functions $g$ and $\psi$ on $\reals^N$ by
\[
g(x) = e^{-\nrm{x}{2}},
\quad
\psi(\xi) = C_N / (1 + 4\pi^2\nrm{\xi}{2}^2)^{(N + 1)/2}
\]
where $C_N$ is the constant $2^N \pi^{(N - 1)/2} \Gamma((N + 1)/2) > 0$.

\begin{lemma}   \label{lemma:sw-ft}
$\psi \in L_1(\reals^N)$ and $\ft{\psi} = g$.
\end{lemma}

\begin{proof}
The first statement is straightforward.  Theorem~1.14 of~\cite{StWe} states
that $\ft{g} = \psi$; but $g$ is continuous and even, so the second statement
follows by Fourier inversion.  \done
\end{proof}

The next lemma is elementary and standard (e.g.~\cite{Katz}).

\begin{lemma}   \label{lemma:easy-Bochner}
Let $\phi \in L_1(\reals^N)$, let $A$ be a finite subset of $\reals^N$, and
let $v \in \reals^A$.  Then
\[
\sum_{a, b \in A} v(a) \ft{\phi}(a - b) v(b)
=
\int_{\reals^N} 
\biggl| \sum_{a \in A} v(a) e^{-2\pi i \ip{\xi}{a}} \biggr|^2 
\phi(\xi)\dee\xi.
\]
\ \done
\end{lemma}

In analytic language, our task is to show that the function $g$ is
strictly positive definite.  This would follow from the easy half of
Bochner's Theorem~\cite{Katz}, except that Bochner's Theorem concerns
\emph{non-strict} positive definiteness.  We therefore need to refine the
argument slightly.

\begin{thm}     \label{thm:l2-pd}
Every finite subspace of Euclidean space is positive definite.
\end{thm}

\begin{proof}
Let $A$ be a finite subspace of $\ell_2^N$.  Let $v \in \reals^A$.  Then
\[
v^\tra \zeta_A v
=
\sum_{a, b \in A} v(a) g(a - b) v(b)
=
\int_{\reals^N} 
\biggl| \sum_{a \in A} v(a) e^{-2\pi i \ip{\xi}{a}} \biggr|^2 
\psi(\xi)\dee\xi
\geq 
0
\]
by Lemmas~\ref{lemma:sw-ft} and~\ref{lemma:easy-Bochner}.
Suppose that $v \neq 0$.  The characters $e^{-2\pi i \ip{\cdot}{a}}$ ($a \in
A$) are linearly independent, so the squared term is positive (that is,
strictly positive) for some $\xi \in \reals^N$.  By continuity, the squared
term is positive for all $\xi$ in some nonempty open subset of $\reals^N$.
Moreover, $\psi$ is continuous and everywhere positive.  So the integral is
positive, as required.  \done
\end{proof}

On the other hand, some of the weights on a finite subspace of Euclidean space
can be negative; see Willerton~\cite{WillHCC} for examples.

\begin{cor}
Every finite subspace of Euclidean space has magnitude.
\done
\end{cor}

A similar argument gives an alternative proof of Theorem~\ref{thm:l1-pd}, that
finite subspaces of $\ell_1^N$ are positive definite.  For this we use the
explicit formula for the Fourier transform of $x \goesto e^{-\nrm{x}{1}}$.
For $p \neq 1, 2$
there is no known formula for the Fourier transform of $e^{-\nrm{x}{p}}$, so
matters become more difficult.  Nevertheless,
Meckes~\cite[Section~\mrk{3}]{MecPDM} has shown that every finite subspace of
$\ell_p^N$ is positive definite whenever $0 < p \leq 2$, and that this is
false for $p > 2$.

\section{Compact metric spaces}
\label{sec:inf}

To extend the notion of magnitude from finite to infinite spaces, there are
broadly speaking two strategies.

In the first, we approximate an infinite space by finite spaces.  As an
initial attempt, given a compact metric space $A$, we might take a sequence
$(A_k)$ of finite metric spaces converging to $A$ in the Gromov--Hausdorff
metric, and try to define $\mg{A}$ as the limit of the sequence $(\mg{A_k})$.
However, this definition is inconsistent; recall Example~\ref{eg:6-over-5}.
We might respond by constraining the sequence $(A_k)$---for example, by taking
$(A_k)$ to be a sequence of subsets of $A$ converging to $A$ in the Hausdorff
metric.

The second strategy is to work directly with the infinite space, replacing
finite sums by integrals.  Weightings are now measures, or perhaps
distributions.  For example, a \demph{weight measure} on a metric space $A$ is
a finite signed Borel measure $w$ such that $\int_A e^{-d(a, b)} dw(b) = 1$
for all $a \in A$.  If $A$ admits a weight measure $w$ then an argument
similar to Lemma~\ref{lemma:wtg-cowtg} shows that $w(A)$ is independent of the
choice of $w$, and we may define the magnitude of $A$ to be $w(A)$.  This was
the definition used by Willerton in~\cite{WillMSS}.

Meckes~\cite{MecPDM} has shown that to a large extent, these different
approaches produce the same result.  Here we implement the first strategy,
defining the magnitude of a space to be the supremum of the magnitudes of its
finite subspaces.  This works well when the space is compact and its finite
subspaces are positive definite.

\subsection{The magnitude of a compact metric space}

\begin{defn}
A metric space is \demph{positive definite} if every finite subspace
is positive definite.  The \demph{magnitude} of a compact positive definite
space $A$ is
\[
\mg{A} 
=
\sup \{ \mg{B} \such B \text{ is a finite subspace of } A \}
\in
[0, \infty].
\]
\end{defn}

These definitions are consistent with the definitions for finite metric
spaces, by Lemma~\ref{lemma:pd-closure}(\ref{part:pd-sub}) and
Corollary~\ref{cor:pd-sub-mag}.

There may even be non-compact spaces for which this definition of magnitude is
sensible.  For example, let $t > 0$, and let $A$ be a space with infinitely
many points and $d(a, b) = t$ for all $a \neq b$; then every finite subspace
of $A$ is positive definite, and the supremum of their magnitudes is $e^t <
\infty$.  In any case, we confine ourselves to compact spaces.

A metric space $A$ is \demph{stably positive definite} if $tA$ is positive
definite for all $t > 0$, or equivalently if every finite subspace of $A$ is
stably positive definite.  (A further equivalent condition, due to
Meckes, is that $A$ is of negative type~\cite[Theorem~\mrk{3.3}]{MecPDM}.)
We already know that $\ell_1^N$ and $\ell_2^N$ are stably positive definite;
much of the rest of this paper concerns the magnitudes of their compact
subspaces.  Ultrametric spaces are also stably positive definite
(Proposition~\ref{propn:um-pd}), and, if compact, have finite magnitude
(Corollary~\ref{cor:um-diam}).  Many other commonly occurring spaces
are stably positive definite too; see~\cite[Theorem~\mrk{3.6}]{MecPDM}.

\begin{defn}
Let $A$ be a stably positive definite compact metric space.  The
\demph{magnitude function} of $A$ is the function
\[
\begin{array}{ccc}
(0, \infty)     &\to            &[0, \infty]    \\
t               &\goesto        &\mg{tA}.
\end{array}
\]
\end{defn}

\begin{lemma}
Let $A$ be a compact positive definite metric space.  Then:
\begin{enumerate}
\item
Every closed subspace $B$ of $A$ is positive definite, and $\mg{B} \leq
\mg{A}$.   
\item
If $A$ is nonempty then $\mg{A} \geq 1$.
\done
\end{enumerate}
\end{lemma}


\begin{propn}   \label{propn:prod-ms-inf}
Let $A$ and $B$ be compact positive definite spaces.  Then $A \otimes B$ is
compact and positive definite, and $\mg{A \otimes B}  = \mg{A} \mg{B}$.
\end{propn}

In the case $A = \emptyset$ and $\mg{B} = \infty$, we interpret $0 \cdot
\infty$ as $0$.

\begin{proof}
Let $C$ be a finite subspace of $A \otimes B$.  Then $C \sub A' \otimes B'$
for some finite subspaces $A' \sub A$ and $B' \sub B$.  Since $A$ and $B$ are
positive definite, so are $A'$ and $B'$.  By Lemma~\ref{lemma:pd-closure}, $A'
\otimes B'$ is positive definite, so $C$ is positive definite.  Hence $A
\otimes B$ is positive definite.
A similar argument shows that $\mg{A \otimes B} = \mg{A} \mg{B}$, using 
Proposition~\ref{propn:prod-ms-fin} and Corollary~\ref{cor:pd-sub-mag}.
\done
\end{proof}

Similarly, Proposition~\ref{propn:union-fin} on unions extends to the compact
setting. 

\begin{propn}   \label{propn:union-inf}
Let $X$ be a metric space and $A, B \sub X$, with $A$ and $B$ compact and $A
\cup B$ positive definite.  Suppose that $A$ projects to $B$ and $B$ projects
to $A$.  Then 
\[
\mg{A \cup B} + \mg{A \cap B}
=
\mg{A} + \mg{B}.
\]
\end{propn}


\begin{proof}
Let $\epsln > 0$.  Choose finite sets
$E \sub A \cup B$ and $H \sub A \cap B$ such that $\mg{A \cup B} \leq \mg{E} +
\epsln$ and $\mg{A \cap B} \leq \mg{H} + \epsln$.  For each $a
\in E \cap A$, choose $\pi_A(a) \in A \cap B$ satisfying the condition of
Definition~\ref{defn:projects-to}, and similarly $\pi_B(b)$ for $b \in E \cap
B$.  Put
\[
H' 
= H \cup \pi_A(E \cap A) \cup \pi_B(E \cap B),
\quad
F  
= (E \cap A) \cup H',
\quad
G = (E \cap B) \cup H'.
\]
Then $F$ and $G$ are finite subsets of $X$, each projecting to the other.
Also $E \sub F \cup G$ and $H \sub F \cap G$.  Applying
Proposition~\ref{propn:union-fin} to $F$ and $G$ gives $\mg{A \cup B} + \mg{A
\cap B} \leq \mg{A} + \mg{B} + 2\epsln$.  Since $\epsln$ was arbitrary, 
%
%
$\mg{A \cup B} + \mg{A \cap B} \leq \mg{A} + \mg{B}$.  

For the opposite inequality, again let $\epsln > 0$, and choose finite sets $F
\sub A$ and $G \sub B$ such that $\mg{A} \leq \mg{F} + \epsln$ and $\mg{B}
\leq \mg{G} + \epsln$.  For each $a \in F$, choose $\pi_A(a) \in A \cap B$
satisfying the condition of Definition~\ref{defn:projects-to}, and similarly
$\pi_B(b)$ for $b \in G$.  Put
\[
F' = F \cup \pi_A F \cup \pi_B G,
\qquad
G' = G \cup \pi_A F \cup \pi_B G.   
\]
Then $F'$ and $G'$ are finite subsets of $X$, each projecting to the other;
also $F \sub F' \sub A$ and $G \sub G' \sub B$.  A similar argument proves
that $\mg{A} + \mg{B} \leq \mg{A \cup B} + \mg{A \cap B} + 2\epsln$.
\done
\end{proof}

\subsection{Subsets of the real line}

As soon as we ask about the magnitude of real intervals, connections with
geometric measure begin to appear.

\begin{propn}   \label{propn:interval-converge}
Let $t \geq 0$ and let $(A_k)$ be a sequence of finite subsets of $\reals$
converging to $[0, t]$ in the Hausdorff metric.  Then $(\mg{A_k})$ converges
to $1 + t/2$.
\end{propn}

This result was announced in~\cite{NCMS}, and also appears, with a different
proof, as Proposition~6 of~\cite{AMSES}.

\begin{proof}
Given $A = \{ a_0 < \cdots < a_n \} \sub \reals$, we have 
\[
(1 + t/2) - \mg{A}
=
\sum_{i = 1}^n
\Bigl\{
\frac{a_i - a_{i - 1}}{2}
-
\tanh\Bigl( \frac{a_i - a_{i - 1}}{2} \Bigr)
\Bigr\}
+
\frac{t - (a_n - a_0)}{2}
\]
by Corollary~\ref{cor:mag-reals}.  The result will follow from the facts that
$\tanh(0) = 0$ and $\tanh'(0) = 1$.  Indeed, write $f(x) = (x - \tanh(x))/x$,
so that $f(x) \to 0$ as $x \to 0$.  Then
\[
\bigl| 
(1 + t/2) - \mg{A}
\bigr|
\leq
\Bigl( \frac{a_n - a_0}{2} \Bigr)
\max_{1 \leq i \leq n} 
\Bigl| f\Bigl( \frac{a_i - a_{i - 1}}{2} \Bigr) \Bigr|
+
\Bigl|
\frac{t - (a_n - a_0)}{2}
\Bigr|.
\]
But $\max_i (a_i - a_{i - 1}) \to 0$
and $a_n - a_0 \to t$ as $A \to [0, t]$, proving the proposition.
\done
\end{proof}

\begin{thm}     \label{thm:mag-interval}
The magnitude of a closed interval $[0, t]$ is $1 + t/2$.
\end{thm}

\begin{proof}
Proposition~\ref{propn:interval-converge} immediately implies that $\mg{[0,
t]} \geq 1 + t/2$.  Now let $A$ be a finite subset of $[0, t]$.  We may choose
a sequence $(A_k)$ of finite subsets of $\reals$ such that $\lim_{k \to
\infty} A_k = [0, t]$ and $A \sub A_k$ for all $k$.  
Then $\mg{A} \leq
\mg{A_k} \to 1 + t/2$ as $k \to \infty$, so
$\mg{A} \leq 1 + t/2$.  
\done
\end{proof}

Ignoring the factor of $1/2$ (which is purely a product of convention),
Theorem~\ref{thm:mag-interval} is a rigorous expression of Schanuel's
contention~\cite{SchaWLP} that the `size' of a closed interval of length
$t$~inches ought to be $(t\text{ inches} + 1)$.

As noted by Willerton~\cite{WillMSS}, there is a weight measure on $[0, t]$.
It is $w = (\delta_0 + \lambda + \delta_t)/2$, where $\delta_x$ is the Dirac
measure at $x$ and $\lambda$ is Lebesgue measure on $[0, t]$.  Then $w([0, t])
= 1 + t/2$.  Hence $w([0, t]) = \mg{[0, t]}$, as guaranteed by
Theorems~\mrk{2.3} and~\mrk{2.4} of Meckes~\cite{MecPDM}.

The magnitude of subsets of $\reals$ is also described by the following
formula, which has no known analogue in higher dimensions.

\begin{propn}
Let $A$ be a compact subspace of $\reals$.  Then
\[
\mg{A} 
=
\frac{1}{2}
\int_\reals
\sech^2 d(x, A) 
\dee x
\]
where $d(x, A) = \inf_{a \in A} d(x, a)$. 
\end{propn}


\begin{proof}
First we prove the identity for finite spaces $A \sub \reals$, by induction on
$n = \card{A}$.  It is elementary when $n \leq 2$.  Now suppose that $n \geq
3$, writing the points of $A$ as $a_1 < \cdots < a_n$.  Put $B = \{a_1,
\ldots, a_{n - 1}\}$ and $C = \{a_{n - 1}, a_n\}$.  Then
\[
\frac{1}{2} \int_\reals \sech^2 d(x, A) \dee x
=
\frac{1}{2} \int_{-\infty}^{a_{n - 1}} \sech^2 d(x, B) \dee x
+
\frac{1}{2} \int_{a_{n - 1}}^\infty \sech^2 d(x, C) \dee x.
\]
Since $\int_0^\infty \sech^2 u \dee u = 1$, this in turn is equal to
\[
\frac{1}{2}
\Bigl( 
\int_\reals \sech^2 d(x, B) \dee x - 1
\Bigr)
+
\frac{1}{2}
\Bigl( 
\int_\reals \sech^2 d(x, C) \dee x - 1
\Bigr)
\]
which by inductive hypothesis is $\mg{B} + \mg{C} - 1$.  On the other hand,
$\mg{A} = \mg{B} + \mg{C} - 1$ by Corollary~\ref{cor:one-point-fin}.  This
completes the induction.  

Now take a compact space $A \sub \reals$.  We know that 
\[
\mg{A} 
=
\sup 
\Bigl\{
\frac{1}{2}
\int_\reals \sech^2 d(x, B) \dee x
\such
B \text{ is a finite subset of } A
\Bigr\}.
\]
Since $\sech^2$ is decreasing on $[0, \infty)$, this implies that 
\[
\mg{A} 
\leq 
\frac{1}{2}\int_\reals \sech^2 d(x, A) \dee x.  
\]
To prove the opposite inequality, choose a sequence $(B_k)$ of finite subsets
of $A$ converging to $A$ in the Hausdorff metric.  We have $0 \leq \sech^2
d(x, B_k) \leq \sech^2 d(x, A)$ for all $x$ and $k$, so 
\[
\lim_{k \to \infty}
\int_\reals
\sech^2 d(x, B_k)
\dee x
=
\int_\reals
\sech^2 d(x, A)
\dee x
\]
by the dominated convergence theorem.  The result follows.
\done
\end{proof}

\subsection{Background on integral geometry}
\label{subsec:int-geom}

To go further, we will need some concepts and results from integral geometry.
Those concerning $\ell_2^N$ can be found in standard texts such
as~\cite{KlRo}.  Those concerning $\ell_1^N$ can be found in~\cite{IGTM}.

Write $\cvx{N}$ for the set of compact convex subsets of $\reals^N$.  A
\demph{valuation} on $\cvx{N}$ is a function $\phi\cln \cvx{N} \to \reals$
such that
%
\begin{equation}        \label{eq:valuation}
\phi(\emptyset) = 0,
\qquad
\phi(A \cup B)
= 
\phi(A) + \phi(B) - \phi(A \cap B)
\end{equation}
whenever $A, B, A \cup B \in \cvx{N}$.  It is \demph{continuous} if continuous
with respect to the Hausdorff metric on $\cvx{N}$, and \demph{invariant} if
$\phi(gA) = \phi(A)$ for all $A \in \cvx{N}$ and isometries $g\cln \ell_2^N
\to \ell_2^N$ (not necessarily fixing the origin).

\begin{examples}
\begin{enumerate}
\item $N$-dimensional Lebesgue measure is a continuous invariant valuation on
$\cvx{N}$, denoted by $\Vol$. 

\item Euler characteristic $\chi$ is a continuous invariant valuation on
$\cvx{N}$.  Since the sets are convex, $\chi(A)$ is $0$ or $1$ according as
$A$ is empty or not.
\end{enumerate}
\end{examples}

The continuous invariant valuations on $\cvx{N}$ form a real vector space,
$\Val{N}$. 

When $A \sub \ell_p^N$ (for any $p \geq 1$) and $t > 0$, the abstract metric
space $tA$ may be interpreted as the subspace $\{ta \such a \in A\}$ of
$\ell_p^N$.  A valuation $\phi$ is \demph{homogeneous} of degree $i$ if
$\phi(tA) = t^i\phi(A)$ for all $A \in \cvx{N}$ and $t > 0$.  

\begin{thm}[Hadwiger~\cite{Hadw}]   \label{thm:Hadwiger}
The vector space $\Val{N}$ has dimension $N + 1$ and a basis $V_0, \ldots,
V_N$ where $V_i$ is homogeneous of degree $i$.
\done
\end{thm}

This description determines the valuations $V_i$ uniquely up to scale factor.
They can be uniquely normalized to satisfy two conditions.  First, $V_N(A) =
\Vol(A)$ for $A \in \cvx{N}$.  Second, whenever $\ell_2^N$ is embedded
isometrically into $\ell_2^{N + 1}$ and $0 \leq i \leq N$, the value $V_i(A)$
is the same whether $A$ is regarded as a subset of $\ell_2^N$ or of $\ell_2^{N
+ 1}$.  With this normalization, $V_i$ is called the $i$th \demph{intrinsic
volume}.

For example, $V_0 = \chi$.  When $A \in \cvx{2}$, $V_1(A)$ is half of the
perimeter of $A$; when $A \in \cvx{3}$, $V_2(A)$ is half of the surface area. 

Here is a general formula for the intrinsic volumes.  For each $0 \leq i \leq
N$, there is an $\Orthog{N}$-invariant measure $\Gram{N}{i}$ on the
Grassmannian $\Gr{N}{i}$, unique up to scale factor.  Given $P \in \Gr{N}{i}$,
write $\pi_P\cln \reals^N \to P$ for orthogonal projection.  Then for $A \in
\cvx{N}$,
\[
V_i(A)
=
c_{N, i}
\int_{\Gr{N}{i}}
\Vol(\pi_P A)\,d\Gram{N}{i}(P)
\]
where $c_{N, i}$ is a positive constant chosen so that the normalizing
conditions are satisfied.

Hadwiger's Theorem solves the classification problem for valuations on
$\ell_2^N$.  More generally, we can try to classify the valuations on any
metric space, in the following sense.

A metric space $A$ is \demph{geodesic}~\cite{Papa} if for all $a, b \in A$
there exists an isometry $ \gamma\cln [0, d(a, b)] \to A $ with $\gamma(0) =
a$ and $\gamma(d(a, b)) = b$.  Given a metric space $X$, write $\cvxa{X}$ for
the set of compact subsets of $X$ that are geodesic with respect to the
subspace metric.  For example, $\cvxa{\ell_2^N} = \cvx{N}$.

A \demph{valuation} on $\cvxa{X}$ is a function $\phi\cln \cvxa{X} \to \reals$
satisfying equations~\eqref{eq:valuation} whenever $A, B, A \cup B, A \cap B
\in \cvxa{X}$.  It is \demph{continuous} if continuous with respect to the
Hausdorff metric, and \demph{invariant} if $\phi(gA) = \phi(A)$ for all
isometries $g$ of $X$.  Write $\Vala{X}$ for the vector space of continuous
invariant valuations on $\cvxa{X}$.  For example, $\Vala{\ell_2^N} = \Val{N}$.

Given any metric space $X$, one can attempt to describe the vector space
$\Vala{X}$.  Here we will need to know the answer for $\ell_1^N$, as well as
$\ell_2^N$.  To state it, we write $\cvxone{N} = \cvxa{\ell_1^N}$ and call its
elements compact \demph{$\ell_1$-convex} sets; similarly, we write $\Valone{N}
= \Vala{\ell_1^N}$.

There are far more $\ell_1$-convex sets than convex sets.  On the other hand,
there are far fewer isometries of $\ell_1^N$ than of $\ell_2^N$; they are
generated by translations, coordinate permutations, and reflections 
in coordinate
hyperplanes.  The following Hadwiger-type theorem is proved
in~\cite{IGTM}.

\begin{thm}
The vector space $\Valone{N}$ has dimension $N + 1$ and a basis $\Vone_0,
\ldots, \Vone_N$ where $\Vone_i$ is homogeneous of degree $i$.
\done
\end{thm}

Again, this determines the valuations $\Vone_i$ uniquely up to scaling.  They
can be described as follows.  For $0 \leq i \leq N$, let $\Grone{N}{i}$ be the
set of $i$-dimensional vector subspaces of $\reals^N$ spanned by some subset
of the standard basis.  For $A \in \cvxone{N}$, put
\[
\Vone_i(A)
=
\sum_{P \in \Grone{N}{i}} \Vol(\pi_P A).
\]
These valuations $\Vone_0, \ldots, \Vone_N$, called the
\demph{$\ell_1$-intrinsic volumes}, satisfy two normalization conditions
analogous to those in the Euclidean case.

The intrinsic volumes of a product space are given by the following formula,
proved in~\cite[Proposition~8.1]{IGTM} and precisely analogous to the
classical Euclidean formula~\cite[Theorem~9.7.1]{KlRo}.

\begin{propn}   \label{propn:intrinsic-prod}
Let $A \in \cvxone{M}$ and $B \in \cvxone{N}$.  Then $A \times B \in \cvxone{M
+ N}$, and
\[
\Vone_k(A \times B)
=
\sum_{i + j = k} \Vone_i(A) \Vone_j(B)
\]
whenever $0 \leq k \leq M + N$.  
\done
\end{propn}

\subsection{Subsets of $\ell_1^N$}

Our investigation of the magnitude of subsets of $\ell_1^N$ begins with sets
of a particularly amenable type.

\begin{defn}
A \demph{cuboid} in $\ell_1^N$ is a subspace of the form $[x_1,
y_1] \times \cdots \times [x_N, y_N]$, where $x_r, y_r \in \reals$ with $x_r
\leq y_r$.
\end{defn}

As an abstract metric space, a cuboid is a tensor product $[x_1, y_1] \otimes
\cdots \otimes [x_N, y_N]$.  

\begin{thm}     \label{thm:mag-cuboid}
For cuboids $A \sub \ell_1^N$, 
\begin{equation}        \label{eq:l1-conj}
\mg{A} 
= 
\sum_{i = 0}^N 2^{-i} \Vone_i(A).
\end{equation}
\end{thm}

\begin{proof}
First let $I = [x, y] \sub \reals$ be a nonempty interval.  By
Theorem~\ref{thm:mag-interval},
\[
\mg{I} 
=
1 + (y - x)/2
=
\chi(I) + \Vol(I)/2
=
\Vone_0(I) + 2^{-1} \Vone_1(I).
\]
This proves the theorem for $N = 1$.  The theorem also holds for $N = 0$.

It now suffices to show that if $A \in \cvxone{M}$ and $B \in \cvxone{N}$
satisfy~\eqref{eq:l1-conj} then so does $A \times B \in \cvxone{M + N}$.
Indeed, as a metric space, $A \times B \sub \ell_1^{M +
N}$ is $A \otimes B$, and the result follows from 
Propositions~\ref{propn:prod-ms-inf} and~\ref{propn:intrinsic-prod}.
\done
\end{proof}

In fact, $\Vone_i(\prod[x_r, y_r])$ is the $i$th elementary symmetric
polynomial in $(y_r - x_r)_{r = 1}^N$, again by
Proposition~\ref{propn:intrinsic-prod}.  It is also equal to
$V_i(\prod[x_r, y_r])$, the Euclidean intrinsic volume.  But in general, the
Euclidean and $\ell_1$-intrinsic volumes of a convex set are not equal.

\begin{cor}     \label{cor:cuboid-mf}
The magnitude function of a cuboid $A \sub \ell_1^N$ is given by
\[
\mg{tA}
=
\sum_{i = 0}^N 2^{-i} \Vone_i(A) t^i.
\]
In particular, the magnitude function of a cuboid $A$ is a polynomial whose
degree is the dimension of $A$, and whose coefficients are proportional to the
$\ell_1$-intrinsic volumes of $A$. 
\done
\end{cor}

The moral is that for spaces belonging to this small class, the dimension and
all of the $\ell_1$-intrinsic volumes can be recovered from the magnitude
function.  In this sense, magnitude subsumes those invariants.  For the rest
of this work we advance the conjectural principle---first set out
in~\cite{AMSES}---that the same is true for a much larger class of spaces, in
both $\ell_1^N$ and $\ell_2^N$.

We begin by showing that the principle holds for subspaces of $\ell_1^N$ when
the invariant concerned is dimension.

\begin{defn}
The \demph{growth} of a function $f\cln (0, \infty) \to \reals$
is 
\[
\inf
\{ \nu \in \reals \such f(t)/t^\nu \text{ is bounded for } t \gg 0 \}
\in 
[-\infty, \infty].
\]  
\end{defn}

For example, the growth of a polynomial is its degree.

\begin{defn}
The \demph{(magnitude) dimension} $\mdim A$ of a stably positive definite
compact metric space $A$ is the growth of its magnitude function.
\end{defn}

\begin{examples}        \label{egs:dim-l1}
\begin{enumerate}
\item   \label{eg:dim-cuboid}
The magnitude dimension of a cuboid in $\ell_1^N$ is its dimension in the
usual sense, by Corollary~\ref{cor:cuboid-mf}.

\item   \label{eg:dim-fin-space}
The magnitude dimension of a nonempty finite space is $0$, by
Proposition~\ref{propn:mf-basic}(\ref{part:mf-lim}). 
\end{enumerate}
\end{examples}

\begin{lemma}   \label{lemma:dim-basic}
Let $A$ be a compact stably positive definite space.  Then:
\begin{enumerate}
\item   \label{part:dim-ineq}
Every closed subspace $B \sub A$ satisfies $\mdim B \leq \mdim A$.
\item   \label{part:dim-nonneg}
If $A \neq \emptyset$ then $\mdim A \geq 0$.
\end{enumerate}
\end{lemma}

\begin{proof}
For~(\ref{part:dim-ineq}), we have $0 \leq \mg{tB} \leq \mg{tA}$ for all $t >
0$, so $\mdim B \leq \mdim A$.  For~(\ref{part:dim-nonneg}), take $B$ to be a
one-point subspace of $A$.  \done
\end{proof}

Recall that the magnitude of a compact positive definite space can in
principle be infinite (although there are no known examples).


\begin{thm}     \label{thm:dim-l1}
Let $A$ be a compact subset of $\ell_1^N$.  Then: 
\begin{enumerate}
\item   \label{part:dim-l1-finite}
$\mg{A} < \infty$.
\item   \label{part:dim-l1-ineq}
$\mdim A \leq N$, with equality if $A$ has nonempty interior.
\end{enumerate}
\end{thm}

We will show in Theorem~\ref{thm:dim-eq} that the hypothesis `nonempty
interior' can be relaxed to `positive measure'.

\begin{proof}
$A$ is a subset of some cuboid $B \sub \ell_1^N$, which has finite
magnitude by Theorem~\ref{thm:mag-cuboid}, so $\mg{A} \leq \mg{B} < \infty$.
Also $\mdim A \leq \mdim B \leq N$ by Lemma~\ref{lemma:dim-basic} and 
Example~\ref{egs:dim-l1}(\ref{eg:dim-cuboid}).  If $A$ has nonempty interior
then it contains an $N$-dimensional cuboid, giving $\mdim A \geq N$.
\done
\end{proof}

We now ask whether the $\ell_1$-intrinsic volumes of an $\ell_1$-convex set
can be extracted from its magnitude function.

Let $\cuboidclass{N}$ be the smallest class of compact subsets of $\ell_1^N$
containing all cuboids and closed under unions of the type in
Proposition~\ref{propn:union-inf}.  By that proposition and
Theorem~\ref{thm:mag-cuboid}, equation~\eqref{eq:l1-conj} holds for all $A \in
\cuboidclass{N}$.

\begin{example}
Let $T$ be a compact triangle in $\ell_1^2$ with two edges
parallel to the coordinate axes (Fig.~\ref{fig:triangle}).
\begin{figure}
\centering
\setlength{\unitlength}{.8em}
\begin{picture}(23.5,9)
\put(0,1.5){\includegraphics[width=23.5\unitlength]{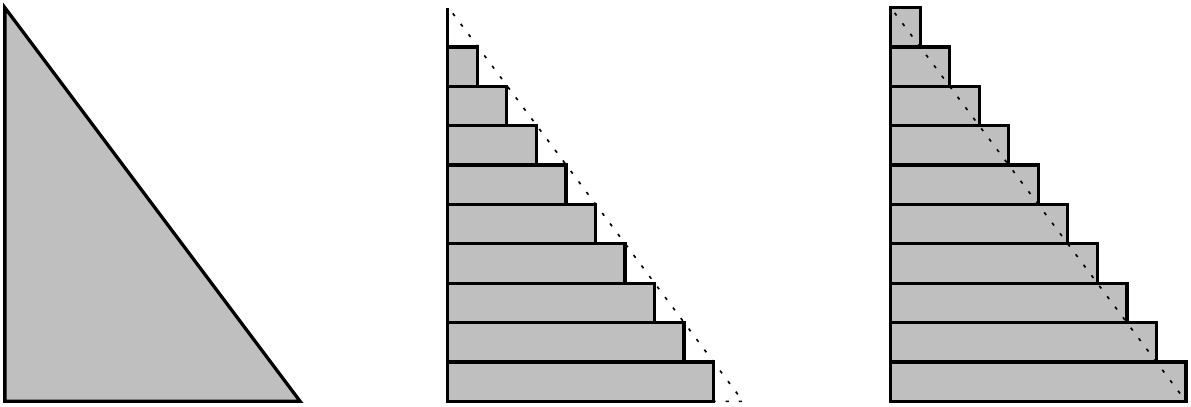}}
\cell{3}{0.15}{b}{T}
\cell{12}{0}{b}{I_{10}}
\cell{21}{0}{b}{E_{10}}
\end{picture}%
\caption{Triangle $T \sub \ell_1^2$, with interior and exterior approximations
$I_{10}$, $E_{10}$}
\label{fig:triangle}
\end{figure}
We compute $\mg{T}$ by exhaustion.  
For each $k \geq 1$, let $I_k$ be the
union of $k$ rectangles approximating $T$ from the inside as in
Fig.~\ref{fig:triangle}; similarly, let $E_k$ be the exterior approximation by
$k$ rectangles.  Then $T$, $I_k$ and $E_k$ are all $\ell_1$-convex with $I_k,
E_k \in \cuboidclass{2}$, and
$\lim_{k \to \infty} I_k = T = \lim_{k \to \infty} E_k$, so
\[
\lim_{k \to \infty} \mg{I_k}
= 
\lim_{k \to \infty} \sum_{i = 0}^2 2^{-i} \Vone_i(I_k)
=
\sum_{i = 0}^2 2^{-i} \Vone_i(T)
=
\lim_{k \to \infty} \sum_{i = 0}^2 2^{-i} \Vone_i(E_k)
=
\lim_{k \to \infty} \mg{E_k}.
\]
But $\mg{I_k} \leq \mg{T} \leq \mg{E_k}$ for all $k$, so
$
\mg{T} = \sum_{i = 0}^2 2^{-i} \Vone_i(T).
$
Similar arguments prove this identity for all compact convex polygons in
$\ell_1^2$.
\end{example}

These and other examples suggest the following conjecture.

\begin{conj}    \label{conj:l1-convex}
Let $A$ be a compact $\ell_1$-convex subspace of $\ell_1^N$.  Then
\[
\mg{A} = \sum_{i = 0}^N 2^{-i} \Vone_i(A).
\]
\end{conj}

If the conjecture holds then $\mg{tA} = \sum_{i = 0}^N 2^{-i} \Vone_i(A)
t^i$ for all $t > 0$ and $A \in \cvxone{N}$.  Hence we can recover all of the
$\ell_1$-intrinsic volumes of an $\ell_1$-convex set from its magnitude
function.

\subsection{Subsets of Euclidean space}

We now prove results for $\ell_2^N$ similar to some of those for
$\ell_1^N$.  
%
Our first task is to prove that the magnitude of a compact subset of Euclidean
space is finite.  Given $A \sub \reals^N$, write
\[
\Sch(A)
=
\{ 
\text{Schwartz functions } \phi\cln \reals^N \to \reals
\text{ such that } \ft{\phi}(a - b) = 1 
\text{ for all } a, b \in A
\}.
\]

\begin{lemma}   \label{lemma:Sch-nonempty}
Let $A \sub \reals^N$ be a bounded set.  Then $\Sch(A) \neq \emptyset$. 
\end{lemma}

\begin{proof}
Since $A$ is bounded, there is a real even Schwartz function $f$ such
that $f(a - b) = 1$ for all $a, b \in A$; then there is a unique real Schwartz
function $\phi$ such that $\ft{\phi} = f$.
\done
\end{proof}

The rest of the proof uses the function $\psi$ from
Section~\ref{subsec:fin-euc}.  For a Schwartz function $\phi$ on $\reals^N$,
write
\[
\fnbound{\phi} 
= 
\sup_{\xi \in \reals^N} 
| \phi(\xi)/\psi(\xi) |
< 
\infty.
\]

\begin{lemma}   \label{lemma:ft-upper-bound}
Let $A$ be a compact subspace of $\ell_2^N$ and $\phi \in \Sch(A)$.  Then
$\mg{A} \leq \fnbound{\phi}$.
\end{lemma}

\begin{proof}
Let $B$ be a finite subset of $A$.  Then for all $v \in \reals^B$, using
Lemma~\ref{lemma:easy-Bochner}, 
\begin{align*}
\fnbound{\phi} \cdot v^\tra \zeta_B v &
=
\fnbound{\phi}
\int_{\reals^N} 
\biggl| \sum_{a \in B} v(a) e^{-2\pi i \ip{\xi}{a}} \biggr|^2 
\psi(\xi)\dee\xi        \\
&
\geq
\int_{\reals^N} 
\biggl| \sum_{a \in B} v(a) e^{-2\pi i \ip{\xi}{a}} \biggr|^2 
\phi(\xi)\dee\xi        
=
\sum_{a, b \in B} v(a) \ft{\phi}(a - b) v(b)   
=
\Bigl( \sum_{a \in B} v(a) \Bigr)^2.
\end{align*}
Taking $v$ to be the weighting on $B$ gives $\fnbound{\phi} \geq
\mg{B}$.  
\done
\end{proof}

\begin{propn}
The magnitude of a compact subspace of $\ell_2^N$ is finite.
\done
\end{propn}


We can extract more from the argument.  For a compact set $A \sub \reals^N$,
write
\[
\spacebound{A} 
= 
\inf \{ \fnbound{\phi} \such \phi \in \Sch(A) \}
< 
\infty.
\]
Lemma~\ref{lemma:ft-upper-bound} states that $\mg{A} \leq \spacebound{A}$.

\begin{lemma}
Let $A$ be a compact subset of $\reals^N$ and $t \geq 1$.  Then
$\spacebound{tA} \leq t^N \spacebound{A}$.
\end{lemma}

\begin{proof}
Let $\phi \in \Sch(A)$.  Define $\theta\cln \reals^N \to \reals$ by
$\theta(\xi) = t^N \phi(t\xi)$.  Then $\theta$ is Schwartz, and if $a, b \in
tA$ then $\ft{\theta}(a - b) = \ft{\phi} ((a - b)/t) = 1$.  Hence $\theta \in
\Sch(tA)$.

I now claim that $\fnbound{\theta} \leq t^N \fnbound{\phi}$.  Indeed, using
the fact that $\psi(\xi) \geq \psi(t\xi)$ for all $\xi \in \reals^N$,
\[
\fnbound{\theta}
=
t^N
\sup_{\xi \in \reals^N} 
\biggl|
\frac{\phi(t\xi)}{\psi(\xi)}
\biggr|
\leq
t^N
\sup_{\xi \in \reals^N}
\biggl|
\frac{\phi(t\xi)}{\psi(t\xi)}
\biggr|
=
t^N \fnbound{\phi}.
\]
This proves the claim, and the result follows.
\done  
\end{proof}

\begin{thm}     \label{thm:dim-l2-upper}
A compact subspace of $\ell_2^N$ has dimension at most $N$.
\end{thm}

\begin{proof}
For compact $A \sub \ell_2^N$ and $t \geq 1$, we have $\mg{tA} \leq
\spacebound{tA} \leq t^N \spacebound{A}$; hence $\mdim A \leq N$.  
\done
\end{proof}

The same proof can be adapted to $\ell_1^N$, although we already have a much
more elementary proof (Theorem~\ref{thm:dim-l1}). 

Having bounded magnitude from above, we now bound it from below.  

\begin{thm}     \label{thm:top-ineq}
Let $\nm{\cdot}$ be a norm on $\reals^N$ whose induced metric is positive
definite.  Write $B = \{ x \in \reals^N \such \nm{x} \leq 1\}$.  For a compact
set $A \sub \reals^N$, equipped with the subspace metric,
\[
\mg{A}
\geq 
\frac{\Vol(A)}{N!\Vol(B)}.
\]
\end{thm}

Before proving this, we state some consequences.  Write
$\omega_N$ for the volume of the unit Euclidean $N$-ball.

\begin{cor}     \label{cor:leading-ineqs}
Let $A$ be a compact subset of $\reals^N$. 
\begin{enumerate}
\item   \label{part:li-l2}
If $A$ is given the subspace metric from $\ell_2^N$ then $\mg{A} \geq
\Vol(A)/N!\omega_N$. 

\item   \label{part:li-l1}
If $A$ is given the subspace metric from $\ell_1^N$ then $\mg{A} \geq
2^{-N} \Vol(A)$. 
\end{enumerate}
\end{cor}

\begin{proof}
Part~(\ref{part:li-l2}) is immediate.  Part~(\ref{part:li-l1}) follows from
the fact that the unit ball in $\ell_1^N$ has volume $2^N/N!$, or can be
derived from Lemma~\ref{lemma:exp-int} below.  \mbox{} \done
\end{proof}

Theorems~\ref{thm:dim-l1}, \ref{thm:dim-l2-upper} and~\ref{thm:top-ineq}
together imply:

\begin{thm}     \label{thm:dim-eq}
Let $p \in \{1, 2\}$ and let $A$ be a compact subspace of $\ell_p^N$.  Then
$\mdim A \leq N$, with equality if $A$ has positive Lebesgue measure.
\done
\end{thm}

Generalizations of these theorems have been proved by Meckes, using more
sophisticated methods~\cite[Theorems~\mrk{4.4} and~\mrk{4.5}]{MecPDM}.  In
particular, Theorem~\ref{thm:dim-eq} is extended to $\ell_p^N$ for all $p \in
(0, 2]$.

To prove Theorem~\ref{thm:top-ineq}, we first need a standard
calculation. 

\begin{lemma}   \label{lemma:exp-int}
Let $\nm{\cdot}$ be a norm on $\reals^N$.  Write $B$ for the unit ball.  Then
\[
\int_{\reals^N} e^{-\nm{x}} \dee x
=
N! \Vol(B).
\]
\end{lemma}

\begin{proof}
$
\displaystyle
\int_{\reals^N} e^{-\nm{x}} \dee x
=
\int_{r = 0}^\infty e^{-r} \dee(\Vol(rB))
=
\int_0^\infty e^{-r} N r^{N - 1} \Vol(B) \dee r
=
N! \Vol(B).
$
\ \done
\end{proof}

\begin{prooflike}{Proof of Theorem~\ref{thm:top-ineq}}
We use the result of Meckes~\cite[Theorem~\mrk{2.4}]{MecPDM} that for a
compact positive definite space $A$ and a finite Borel measure $v$ on $A$,
\[
\mg{A} \geq v(A)^2
\Bigl/
\int_A \int_A e^{-d(a, b)} \dee v(a) \dee v(b).
\]
Let $A \sub \reals^N$ be a compact set and take $v$ to be Lebesgue measure:
then
\begin{align*}
\mg{A}  
\geq 
\frac{\Vol(A)^2}{\int_A \int_A e^{-\nm{a - b}} \dee a \dee b}    
\geq 
\frac{\Vol(A)^2}%
{\int_A \int_{\reals^N} e^{-\nm{a - b}} \dee a \dee b}          
= 
\frac{\Vol(A)^2}%
{\int_A \int_{\reals^N} e^{-\nm{x}} \dee x \dee b}
= 
\frac{\Vol(A)}{\int_{\reals^N} e^{-\nm{x}} \dee x}.
\end{align*}
The theorem follows from Lemma~\ref{lemma:exp-int}.
\done
\end{prooflike}

This proof is a rigorous rendition of part of Willerton's bulk approximation
argument~\cite{WillHCC}.  There is an alternative proof in the same spirit,
not depending on the results of Meckes but instead working with
finite approximations.  We sketch it now.

\begin{prooflike}{Alternative proof of Theorem~\ref{thm:top-ineq}}
For $\delta > 0$, write
\[
S_\delta 
=
\Bigl\{ 
x \in \delta\integers^N
\such
A \cap \prod_{r = 1}^N [x_r, x_r + \delta) \neq \emptyset 
\Bigr\}.
\]
Define $\alpha\cln \delta\integers^N \to \reals^N$ by choosing for each $x \in
S_\delta$ an element $\alpha(x) \in A \cap \prod [x_r, x_r + \delta)$, and
putting $\alpha(x) = x$ for $x \in \delta\integers^n \without S_\delta$.

A calculation similar to that in the first proof of Theorem~\ref{thm:top-ineq}
shows that for all $\delta > 0$,
\[
\mg{A} 
\geq
\frac{\card{S_\delta}}{\sum_{x \in \delta\integers^N} E_\delta(x)}
\]
where
\[
E_\delta(x)
=
\frac{1}{\card{S_\delta}} 
\sum_{y \in S_\delta} 
e^{-\nm{\alpha(x + y) - \alpha(y)}}
\quad
\bigl(
\approx
e^{-\nm{x}}
\bigr).
\]
(Apply Proposition~\ref{propn:CS} to the finite space $\alpha
S_\delta$.)  Since Lebesgue measure is outer regular, $\lim_{\delta \to 0}
(\delta^N (\card{S_\delta})) = \Vol(A)$.  From the fact that
$\nm{\alpha(x) - x} \leq \diam([0, \delta)^N)$ for all $\delta > 0$ and $x \in
\delta\integers^N$, it also follows that
\[
\lim_{\delta \to 0}
\Bigl(
\delta^N 
\!\!
\sum_{x \in \delta\integers^N} E_\delta(x)
\Bigr)
=
\int_{\reals^N} e^{-\nm{x}} \dee x.
\]
The theorem now follows from Lemma~\ref{lemma:exp-int}.
\done
\end{prooflike}

These results suggest the following conjecture, first stated in~\cite{AMSES}:

\begin{conj}    \label{conj:convex}
Let $A$ be a compact convex subspace of $\ell_2^N$.  Then
\[
\mg{A} = 
\sum_{i = 0}^N \frac{1}{i!\omega_i} V_i(A).
\] 
\end{conj}
    
Assuming the conjecture, the magnitude function of a compact convex set $A
\sub \ell_2^N$ is a polynomial:
\begin{equation}        \label{eq:l2-conj}
\mg{tA}
=
\sum_{i = 0}^N \frac{1}{i!\omega_i} V_i(A) t^i.
\end{equation}
All of the intrinsic volumes, as well as the dimension, can therefore be
recovered from the magnitude function.  

The evidence for Conjecture~\ref{conj:convex} is as follows.
\begin{itemize}
\setlength{\itemsep}{0ex}
\item The two sides of equation~\eqref{eq:l2-conj} have the
same growth (by Theorem~\ref{thm:dim-eq}).

\item The left-hand side of~\eqref{eq:l2-conj} is greater than or equal to the
leading term of the right-hand side (by Corollary~\ref{cor:leading-ineqs}).

\item The conjecture holds for $N = 1$ (by Theorem~\ref{thm:mag-interval}).


\item It is closely analogous to Conjecture~\ref{conj:l1-convex}, which, while
itself a conjecture, is known to hold for a nontrivial class of examples.  
(To see the analogy, note that in both cases the $i$th coefficient is
$1/i!\Vol(B_i)$, where $B_i$ is the $i$-dimensional unit ball.)

\item There is good numerical evidence, due to Willerton~\cite{WillHCC}, when
$A$ is a disk, square or cube.
\end{itemize}

One strategy for proving Conjecture~\ref{conj:convex} would be to apply
Hadwiger's Theorem~(\ref{thm:Hadwiger}).  There are currently two obstacles.
First, it is not known that magnitude is a valuation on compact convex sets.
Certainly it is not a valuation on \emph{all} compact subsets of $\ell_2^N$:
consider the union of two points.

Second, even supposing that magnitude is a valuation on convex sets, the
conjecture is not proved.  We would know that magnitude was an invariant
valuation, monotone and therefore continuous by Theorem~8 of
McMullen~\cite{McMVET}.%
\footnote{I thank Mark Meckes for this observation.}
By Hadwiger's Theorem, there would be constants $c_i$ such that $\mg{A} = \sum
c_i V_i(A)$ for all convex sets $A$.  However, current techniques provide no
way of computing those constants.  Knowing the magnitude of balls or
cubes would be enough.  But apart from subsets of the line, there is not
a single convex subset of Euclidean space whose magnitude is known.

\ucontents{section}{References}

\end{document}